\documentclass[11pt,a4paper,leqno,twoside, headinclude]{amsart}

\usepackage[colorlinks,pdftex, plainpages=false]{hyperref}
\hypersetup{pdftitle=Partial Classification Results for Positive
Quaternion K\"ahler Manifolds, pdfauthor=Manuel Amann,
pdftoolbar=false, plainpages=false, hyperindex=true,
pdfdisplaydoctitle=true}

\hypersetup{colorlinks,%
citecolor=black,%
filecolor=black,%
linkcolor=black,%
urlcolor=black,%
pdftex}

\usepackage{color}

\usepackage{amsmath, amssymb, amscd, amsthm}
\usepackage{stmaryrd}
\usepackage{latexsym}
\usepackage[all]{xy}
\usepackage{pb-diagram}
\usepackage{rotating}
\usepackage{multicol}
\usepackage{lscape}
\usepackage{wasysym}
\usepackage{longtable}

\xyoption{all}

\newtheorem{theo}{Theorem}[section]
\newtheorem{main}{Theorem}
\newtheorem*{main*}{Theorem}

\newtheorem{prop}[theo]{Proposition}
\newtheorem{defi2}[theo]{Definition}
\newenvironment{defi}{\begin{defi2}\normalfont}{\end{defi2}}
\newtheorem{rem2}[theo]{Remark}
\newenvironment{rem}{\begin{rem2}\normalfont}{\hfill$\boxbox$\end{rem2}}
\newtheorem{lemma}[theo]{Lemma}
\newtheorem{cor}[theo]{Corollary}

\newtheorem*{conj*}{Conjecture}
\newtheorem*{theo*}{Theorem}
\newtheorem*{mi2}{Main Idea}

\newtheorem{ex2}[theo]{Example}

\newtheorem{exer2}[theo]{Exercise}

\newtheorem{alg2}[theo]{Algorithm}

\newcommand{\cc}{{\mathbb{C}}}                                     
\newcommand{\hh}{{\mathbb{H}}}                                     
\newcommand{\oo}{{\mathbb{O}}}                                     
\newcommand{\nn}{{\mathbb{N}}}                                     
\newcommand{\qq}{{\mathbb{Q}}}                                     
\newcommand{\rr}{{\mathbb{R}}}                                     
\newcommand{\pp}{{\mathbf{P}}}                                     
\newcommand{\Gr}{{\mathbf{Gr}}}                                    
\newcommand{\s}{{\mathbb{S}}}                                      
\newcommand{\zz}{{\mathbb{Z}}}                                     
\newcommand{\SO}{{\mathbf{SO}}}                                    
\newcommand{\U}{{\mathbf{U}}}                                      
\newcommand{\SU}{{\mathbf{SU}}}                                    
\newcommand{\Sp}{{\mathbf{Sp}}}                                    
\newcommand{\A}{{\mathbf{A}}}                                      
\newcommand{\B}{{\mathbf{B}}}                                      
\newcommand{\C}{{\mathbf{C}}}                                      
\newcommand{\D}{{\mathbf{D}}}                                      
\newcommand{\E}{{\mathbf{E}}}                                      
\newcommand{\F}{{\mathbf{F}}}                                      
\newcommand{\G}{{\mathbf{G}}}                                      
\newcommand{\Spin}{{\mathbf{Spin}}}                                
\newcommand{\In} {{\,\subseteq\,}}                                 
\newcommand{\ch}{{\operatorname{ch}}}                              
\newcommand{\Aut}{{\operatorname{Aut}}}                            
\newcommand{\Isom}{{\operatorname{Isom}}}                          
\newcommand{\Hol}{{\operatorname{Hol}}}                            
\newcommand{\id}{{\operatorname{id}}}                              
\newcommand{\ind} {{\operatorname{ind}}}                           
\newcommand{\rk}{{\operatorname{rk\,}}}                            
\newcommand{\sign}{{\operatorname{sign}}}                          
\newcommand{\diag}{{\operatorname{diag}}}                          
\newcommand{\Dc}{\mbox{$D\!\!\!\! / $}}                            
\newcommand{\vproof}{{\begin{flushright} \qed                      
                      \end{flushright}}}

\newcommand{\comment}[1]{}                                         
\newcommand{\hto}[1]{\overset{#1}{\hookrightarrow}}                
\newcommand{\case}[1]{\textbf{Case #1.}}                           
\newcommand{\step}[1]{\textbf{Step #1.}}                           
\newcommand{\ack}{\noindent\textbf{Acknowledgements. }}            
\newcommand{\str}{\noindent\textbf{Structure of the article. }}    

\newenvironment{prf}{\begin{proof}[\textsc{Proof}]} {\end{proof}}     

\title[Partial Classification Results for PQK Manifolds]{Partial Classification Results for Positive
Quaternion K\"ahler Manifolds}
\author{Manuel Amann}
\date{October 29th, 2009}

\begin{document}

\maketitle \thispagestyle{empty}


\begin{abstract}
Positive Quaternion K\"ahler Manifolds are Riemannian manifolds with
holonomy contained in $\Sp(n)\Sp(1)$ and with positive scalar
curvature. Conjecturally, they are symmetric spaces. We prove this
conjecture in dimension $20$ under additional assumptions and we
provide recognition theorems for quaternionic projective spaces (in
low dimensions) as well as the real Grassmanian $\widetilde
\Gr_4(\rr^{n+4})$.
\end{abstract}


\section*{Introduction}

Quaternion K\"ahler Manifolds settle in the highly remarkable class
of special geometries. Hereby one refers to Riemannian manifolds
with special holonomy among which K\"ahler manifolds, Calabi--Yau
manifolds or Joyce manifolds are to be mentioned as the most
prominent examples. Quaternion K\"ahler Manifolds have holonomy
contained in $\Sp(n)\Sp(1)$; they are called \emph{positive}, if
their scalar curvature is positive.

The only known examples of Positive Quaternion K\"ahler Manifolds
are given by the so-called \emph{Wolf-spaces}, which are all
symmetric and the only homogeneous examples due to Alekseevski.
Indeed, they are given by the infinite series $\hh\pp^n$,
$\Gr_2(\cc^{n+2})$ and $\widetilde \Gr_4(\rr^{n+4})$ (the
Grassmanian of oriented real $4$-planes) and the exceptional spaces
$\G_2/\SO(4)$, $\F_4/\Sp(3)\Sp(1)$, $\E_6/\SU(6)\Sp(1)$,
$\E_7/\Spin(12)\Sp(1)$, $\E_8/\E_7\Sp(1)$. Besides, it is known that
in each dimension there are only finitely many Positive Quaternion
K\"ahler Manifolds. This endorses
the fundamental conjecture
\begin{conj*}[LeBrun, Salamon]\label{PQKconj01}
Every Positive Quaternion K\"ahler Manifold is a Wolf space.
\end{conj*}
A confirmation of the conjecture has been achieved in dimensions
four (Hitchin) and eight (Poon--Salamon, LeBrun--Salamon).

Recently this field of study has received a lot of attention with
several contributions via completely different approaches and
methods ranging from Ricci-flow to complex geometry. However, the
LeBrun--Salamon conjecture still seems to be open. This article is
devoted to an investigation of low-dimensional Positive Quaternion
K\"ahler Manifolds---dimensions $16$ to $24$, with a clear emphasis
on dimension $20$---as well as high-dimensional ones with
symmetries.

In this article we provide several classification results under
different assumptions coming from Index Theory, algebraic topology
or the study of symmetries.

\vspace{3mm}

There used to be a classification of Positive Quaternion K\"ahler
Manifolds in dimension $12$ by Hayde\'e and Rafael Herrera
(cf.~\cite{herr}) confirming the main conjecture. In the second part
of that article they showed that any $12$-dimensional Positive
Quaternion K\"ahler Manifold $M$ is symmetric if the $\hat A$-genus
of $M$ vanishes. If $M$ is a spin manifold, this condition is always
fulfilled by a classical result of Lichnerowicz, since a Positive
Quaternion K\"ahler Manifold has positive scalar curvature. Positive
Quaternion K\"ahler Manifolds $M^n\not\cong \hh\pp^n$ are known to
admit a spin structure if and only if $n$ is even (cf.~proposition
\cite{sal1}.2.3, p.~148). One also knows that $\hat A(M)[M]$
vanishes on the symmetric examples with finite second homotopy group
(cf.~theorem \cite{BoHiII}, 23.3). Atiyah and Hirzebruch
(cf.~\cite{AtHi}) showed that the $\hat A$-genus vanishes on spin
manifolds with smooth effective $\s^1$-action.

In the first part of their article \cite{herr} Hayde\'e and Rafael
Herrera claim a similar result for simply-connected manifolds with
finite second homotopy group instead of a spin structure.
Unfortunately, we found the proof of this assertion to be erroneous
and---as we show together with Anand Dessai in \cite{mam1}---it is
not true in general that the $\hat A$-genus of such a
simply-connected $\pi_2$-finite manifold with smooth effective
$\s^1$-action vanishes. Thus the classification in dimension twelve
can no longer be sustained---see also \cite{herrerr}. However, it
still remains an open question whether the $\hat A$-genus vanishes
on $\pi_2$-finite Positive Quaternion K\"ahler Manifolds.

These observations actually where a byproduct of our work in
dimension $20$---with rather detrimental effects on our own
arguments and theorems. This is the reason why the case of dimension
$12$ is neglected in this article and why the assumption $\hat
A(M)[M]=0$ appears in
\begin{main}\label{theoA}
A $20$-dimensional Positive Quaternion K\"ahler Manifold $M$ with
$\hat A(M)[M]=0$ satisfying
\begin{align*}
\dim \Isom(M)\not\in \{15,22,29\}
\end{align*}
(where possible groups $\Isom_0(M)$ in these dimensions can be read
off from table \ref{INTtable01}) is a Wolf space.
\end{main}
The proof of this theorem essentially splits into two parts: On the
one hand we combine relations from Index Theory with further
properties of Positive Quaternion K\"ahler Manifolds to restrict
possible identity components of the isometry group to a small list
of relatively large groups. In the second part we use Lie theoretic
arguments to provide the classification result. Both approaches path
the way towards further results:

\vspace{3mm}

In the vein of the index computations we obtain a classification
result concerning the quaternionic projective spaces.
\begin{main}\label{theoB}
Suppose $2 \neq n\leq 6$. Let $M^{4n}$ be a Positive Quaternion K\"ahler
Manifold with $b_4(M)=1$. Then $M$ is homothetic to $\hh\pp^n$.
\end{main}
This result has been proven in dimensions $\dim M\leq 16$ in
\cite{sas}. The methods applied there permit a generalisation to
dimensions $20$ and $24$. Note that the exceptional Wolf space
$\F_4/\Sp(3)\Sp(1)$ has dimensions $28$ (and $\dim
\G_2/\SO(4)=8$)---both of them satisfy $b_4=1$.

The arguments arising from the theory of transformation groups can
be generalised to yield a recognition theorem for the real
Grassmanian by means of the dimension of its isometry group. The
complex Grassmanian and the quaternionic projective space are
topologically well-identifiable by means of the second homotopy
group, i.e.~$\pi_2(M^{4n})=0$ implies $M^{4n}\cong \hh\pp^n$ and
$\pi_2(M)=\zz$ leads to $M^{4n}\cong \Gr_2(\cc^{n+2})$. These
identities make it possible to recognize these two spaces by means
of the rank of the isometry group (cf.~\cite{fang2}). No similar
characterisations seem to be known for the real Grassmanian which
makes our recognition theorem a first one of its kind.
\begin{main}\label{theoC}
Let $M^{4n}$ be a Positive Quaternion K\"ahler Manifold. Suppose
that the dimension of the isometry group $\dim \Isom(M^{4n})$
satisfies the respective condition depicted in table
\ref{tableC}. Then $M$ is symmetric and it holds:
\begin{align*}
\begin{array}{ll}
M\cong \widetilde \Gr_4(\rr^{n+4}) & \Leftrightarrow \dim
\Isom(M)=\frac{n^2+7n+12}{2} \\&\Leftrightarrow
\rk \Isom(M)=\big\lfloor \frac{n}{2} \big\rfloor+2\\
M\cong \Gr_2(\cc^{n+2}) &\Leftrightarrow \dim \Isom(M)=n^2+4n+3\\
M\cong \hh\pp^n &\Leftrightarrow \dim \Isom(M)=2n^2+5n+3
\end{array}
\end{align*}
In particular, if  the isometry group satisfies that
\begin{align*}
\dim \Isom(M)> \frac{n^2+5n+12}{2}
\end{align*}
for $n\geq 22$ and $n\not\in \{27,28\}$, then $M$ is symmetric and
we recognise the real Grassmannian by the dimension of its isometry
group.
\end{main}
\begin{table}[h]
\centering \caption{A recognition theorem} \label{tableC}
\begin{tabular}{|@{\hspace{6mm}}r@{\hspace{6mm}}| @{\hspace{6mm}}r@{\hspace{6mm}}|@{\hspace{6mm}}c@{\hspace{6mm}}| }
\hline $n=$ & $\dim \Isom(M^{4n})>$ & recognising  \\
\hline\hline $3$ & $28$ & $\hh$ \\
\hline $4$ & $52$ & $\hh$ \\
\hline $5$ & $55$ & $\hh$ \\
\hline $6$ & $55$ & $\hh,\cc$ \\
\hline $7$ & $78$ & $\hh,\cc$ \\
\hline $8$ & $78$ & $\hh,\cc$ \\
\hline $9$ & $133$ & $\hh$ \\
\hline $10$ & $133$ & $\hh,\cc$ \\
\hline $11$ & $248$ & $\hh$ \\
\hline $12$ & $248$ & $\hh$ \\
\hline $13$ & $251$ & $\hh$ \\
\hline $14$ & $251$ & $\hh,\cc$ \\
\hline $15$ & $262$ & $\hh,\cc$ \\
\hline $16$ & $262$ & $\hh,\cc$ \\
\hline $17$ & $269$ & $\hh,\cc$ \\
\hline $18$ & $269$ & $\hh,\cc$ \\
\hline $19$ & $300$ & $\hh,\cc$ \\
\hline $20$ & $300$ & $\hh,\cc$ \\
\hline $21$ & $303$ & $\hh,\cc$ \\
\hline
\hline $22$ & $303$ & $\hh,\cc,\rr$ \\
\hline $23$ & $328$ & $\hh,\cc,\rr$ \\
\hline $24$ & $354$ & $\hh,\cc,\rr$ \\
\hline $25$ & $381$ & $\hh,\cc,\rr$ \\
\hline $26$ & $409$ & $\hh,\cc,\rr$ \\
\hline $27$ & $496$ & $\hh,\cc$ \\
\hline $28$ & $496$ & $\hh,\cc$ \\
\hline
\hline $n\geq 29$ & $\frac{n^2+5n+12}{2}$ & $\hh,\cc,\rr$ \\
\hline
\end{tabular}
\end{table}
The symbols $\hh,\cc,\rr$ in the column ``recognising'' in table
\ref{tableC} refer to whether we may identify the quaternionic
projective space, the complex Grassmannian or the real Grassmannian
in this dimension by the theorem.

We remark that in general the dimension of the isometry group
\linebreak[4]$\dim \Isom(M^{4n})$ takes values in
\begin{align*}
[0,\dim \Sp(n+1)]=[0,2n^2+5n+3]
\end{align*}
(cf.~\ref{PQKtheo07}). Hence---apart from recognising the real
Grassmannian---this theorems rules out approximately three quarters
of all possible values.

As the dimension of the isometry group of a Positive Quaternion
K\"ahler Manifold may be interpreted as the index of a certain
twisted Dirac operator (cf.~theorem \ref{PQKtheo10}) we can make the
following observation. \vspace{3mm}
\begin{center}
{The question whether a Positive Quaternion K\"ahler Manifold
$M^{4n}$ is symmetric or not can (almost always) be decided from the
index}
\begin{align*}
\ind (\Dc(S^{n+2}H)) = \langle\hat A(M) \cdot
\ch(S^{n+2}H),[M]\rangle
\end{align*}
\end{center}
(For the bundle $H$ we refer to the next section.)

\vspace{3mm}

\str In section \ref{sec01} we shall give a very brief introduction
to Positive Quaternion K\"ahler Geometry focussing on properties
obtained via Index Theory or transformation groups. Section
\ref{sec02} is devoted to computations of several twisted $\hat
A$-genera via characteristic classes. In section \ref{sec03} we
shall prove theorem \ref{theoB}. Moreover, we shall identify further
properties related to the methods of the proof. Section \ref{sec04}
will be used to present further results that come out of our index
computations. On the one hand we present theorems (e.g.~on the
existence of isometric $\s^1$-actions) that are of interest for
their own sake. On the other hand this section will establish the
first part of the proof of theorem \ref{theoA}; namely it will yield
the existence of large isometry groups under mild assumptions.
Section \ref{sec05} is devoted to classification of $20$-dimensional
Positive Quaternion K\"ahler Manifolds with large isometry groups.
The results here combine with the ones from chapter \ref{sec04} to
complete the proof of theorem \ref{theoA}. The proof of theorem
\ref{theoC} will be given in \ref{sec06}.

\vspace{3mm}

We remark that a more elaborate introduction to the subject as well
as detailed proofs can be found in \cite{phd}.

Several arguments involve heavy computations. All of these were done
with the help of \textsc{Mathematica 6.01} or \textsc{Maple 9} or
later programmme versions respectively.

\vspace{3mm}

\ack The author is very grateful to Anand Dessai for various
fruitful discussions and suggestions.


\section{Positive Quaternion K\"ahler Manifolds}\label{sec01}

Due to Berger's celebrated theorem the holonomy group $\Hol(M,g)$ of
a simply-connected, irreducible and non-symmetric Riemannian
manifold $(M,g)$ is one of $\SO(n)$, $\U(n)$, $\SU(n)$, $\Sp(n)$,
$\Sp(n)\Sp(1)$, $\G_2$ and $\Spin(7)$.

A connected oriented Riemannian manifold $(M^{4n},g)$ is called a
\emph{Quaternion K\"ahler Manifold} if
\begin{align*}
\Hol(M,g)\In \Sp(n)\Sp(1)=\Sp(n)\times \Sp(1)/\langle -\id,
-1\rangle
\end{align*}
(In the case $n=1$ one additionally requires $M$ to be Einstein and
self-dual.) Quaternion K\"ahler Manifolds are Einstein
(cf.~\cite{bess}.14.39, p.~403). In particular, their scalar
curvature is constant.
\begin{defi}\label{PQKdef01}
A \emph{Positive Quaternion K\"ahler Manifold} is a Quaternion
K\"ahler Manifold with complete metric and with positive scalar
curvature.
\end{defi}
For an elaborate depiction of the subject we recommend the survey
articles \cite{sal1} and \cite{sal2}. We shall content ourselves
with mentioning a few properties that will be of importance
throughout the article:

Foremost, we note that Positive Quaternion K\"ahler Manifolds $M$
clearly are not necessarily K\"ahlerian, as the name might suggest.
Moreover, the manifold $M$ is compact and simply-connected
(cf.~\cite{sal1}, p.~158 and \cite{sal1}.6.6, p.~163).

\vspace{3mm}

Locally the structure bundle with fibre $\Sp(n)\Sp(1)$ may be lifted
to its double covering with fibre $\Sp(n)\times \Sp(1)$. The bundles
associated to the standard complex representations of $\Sp(n)$ on
$\cc^{2n}$ and of $\Sp(1)$ on $\cc^2$ will be called $E$
respectively $H$. Recall that $\pp_\cc(H)$ is called the
\emph{twistor space} of $M$. This space is a Fano contact K\"ahler
Einstein manifold (cf.~theorem \cite{brunsal}.1.2, p.~113).

We obtain the following formula for the complexified
tangent bundle $T_\cc M$ of the Positive Quaternion K\"ahler
Manifold $M$ (cf.~\cite{sal2}, p.~93):
\begin{align*}
T_\cc M=E\otimes H
\end{align*}
The bundles $E$ and $H$ arise from self-dual representations and so
their odd-degree Chern classes vanish. The Chern classes of $E$ will
be denoted by
\begin{align*}
c_{2i}&:=c_{2i}(E)\in H^{2i}(M)\intertext{and} u&:=-c_2(H)\in H^4(M)
\end{align*}
The \emph{quaternionic volume}
\begin{align*}
v=(4u)^n\in H^{4n}(M^{4n})
\end{align*}
is integral and satisfies $1\leq v\leq 4^n$---cf.~\cite{sal2},
p.~114 and corollary \cite{sem1}.3.5, p.~7.

Recall the definiteness of the intersection form on Positive
Quaternion K\"ahler Manifolds (cf.~\cite{fuj}), \cite{naga}) (which
is a consequence of the Hodge--Riemann bilinear relations on the
twistor space). The quoted articles differ in the formulation of
positive/negative definiteness. We use the form $u=-c_2(H)$ to state
the theorem. The orientation on $M$ is naturally given by $u^n$.
\begin{theo}\label{PQKtheo08}
The generalised intersection form
\begin{align*}
Q(x,y)=(-1)^{r/2}\int_M x\wedge y\wedge u^{n-r/2}
\end{align*}
for $[x],[y]\in H^{r}(M^{4n},\rr)$ with even $r\geq 0$ is positive
definite. In particular, the signature of the manifold satisfies
\begin{align*}
\sign(M)=(-1)^n b_{2n}(M)
\end{align*}
\end{theo}
\vproof

An oriented compact manifold $M^{4n}$ is called \emph{spin} if its
$\SO(4n)$-structure bundle lifts to a $\Spin(4n)$-bundle, or
equivalently, its second Stiefel-Whitney class $w_2(M)=0$ vanishes.
A vector bundle $E\to M$ is called \emph{spin} if $w_1(E)=w_2(E)=0$.
\label{page02}Positive Quaternion K\"ahler Manifolds $M^{4n}\neq
\hh\pp^n$ are spin if and only if $n$ is even (cf.~proposition
\cite{sal1}.2.3, p.~148).

On a Positive Quaternion K\"ahler Manifold $M^{4n}$ we have the
locally associated bundles $E$ and $H$ from above. Now form the
(virtual) bundles
\begin{align*}
\sideset{}{_0^{\!\!k}}\bigwedge E&:=\sideset{}{\!\!^k}\bigwedge E -
\sideset{}{\!\!^{k-2}}\bigwedge E \intertext{of exterior powers and
the bundles} S^l H&:=\operatorname{Sym}^l H
\end{align*}
of symmetric powers. In general, the bundles
$\sideset{}{_0^k}\bigwedge E\otimes S^l H$ exist globally if and
only if $n+k+l$ is even. In this case, using the index theorem one
obtains the following relations (cf.~\cite{sal2}, p.~117) where
$i^{k,l}:=\ind \Dc\big(\sideset{}{_0^k}\bigwedge E\otimes S^l
H\big)$ is the index of a twisted Dirac operator.
\begin{theo}\label{PQKtheo10}
It holds:
\begin{align*}
i^{k,l}=
\begin{cases}
0 & \textrm{if } k+l<n\\
(-1)^k (b_{2k}(M)+b_{2k-2}(M)) &\textrm{if } k+l=n \\
d &\textrm{if } k=0,~ l=n+2
\end{cases}
\end{align*}
where $d=\dim \Isom(M)$ is the dimension of the isometry group of
$M$ and the $b_i(M)$ are the Betti numbers of $M$ as usual.
\end{theo}
\vproof
As a consequence of the \emph{Atiyah--Singer Index Theorem}
we then may express these indices topologically via genera:
\begin{align}\label{PQKeqn06}
\ind \bigg(\Dc\bigg(\sideset{}{_0^k}\bigwedge E\otimes S^l
H\bigg)\bigg) = \bigg\langle\hat A(M) \cdot
\ch\bigg(\sideset{}{_0^k}\bigwedge E\bigg)\cdot\ch\bigg(S^l
H\bigg)\bigg),[M]\bigg\rangle
\end{align}

\vspace{3mm}

Let us now collect some information on isometry groups:
\begin{theo}\label{PQKtheo07}
Let $M^{4n}$ be a Positive Quaternion K\"ahler Manifold with
isometry group $\Isom(M)$. We obtain:
\begin{itemize}
\item
The rank $\rk\Isom(M)$ may not exceed $n+1$. If $\rk \Isom(M)=n+1$,
then $M\in\{\hh\pp^n,\Gr_2(\cc^{n+2})\}$.
\item
If $\rk \Isom(M)\geq \frac{n}{2}+3$, then $M$ is isometric to
$\hh\pp^n$ or to $\Gr_2(\cc^{n+2})$.
\item
It holds that $\dim \Isom(M^{4n})\leq \dim \Sp(n+1)=(n+1)(2n+3)$.
Equality holds if and only if $M\cong \hh\pp^n$.
\item
If $n=3$, then $\dim \Isom(M)\geq 5$; if $n=4$, then $\dim
\Isom(M)\geq 8$.
\end{itemize}
\end{theo}
\begin{prf}
The first assertion is due to theorem \cite{sal2}.2.1, p.~89. The
second item is theorem \cite{fang2}.1.1, p.~642. The inequality in
the third assertion follows from corollary \cite{sem1}.3.3, p.~6. In
case $\dim \Isom(M)=(n+1)(2n+3)$ it was already observed on
\cite{sal1}, p.~161 that $M$ is homothetic to the quaternionic
projective space. The fourth point is due to theorem
\cite{sal1}.7.5, p.~169.
\end{prf}
The isometry group $\Isom(M)$ of $M$ is a compact Lie group. Due to
theorems \cite{brock}.V.8.1, p.~233, and \cite{brock}.V.7.13,
p.~229, we may assume up to finite coverings that $\Isom_0(M)$---the
component of the identity---is the product of a simply-connected
semi-simple Lie group and a torus.

For the convenience of the reader we give a table of simple Lie
groups by table \ref{PQKtable04}, which will support our future
arguments involving dimensions and ranks of Lie groups.
\begin{table}[h]
\centering \caption{Classification of simple Lie groups up to
coverings} \label{PQKtable04}
\begin{tabular}{|l@{\hspace{3mm}}| @{\hspace{3mm}}l@{\hspace{3mm}}|
@{\hspace{3mm}}r@{\hspace{3mm}}|}
\hline type & corresponding Lie group & dimension \\
& (not necessarily simply-connected) &\\
\hline\hline $\A_n$, $n=1,2,\dots$ & $\SU(n+1)$ & $n(n+2)$\\
\hline $\B_n$, $n=1,2,\dots$ & $\SO(2n+1)$ & $n(2n+1)$
\\
\hline
$\C_n$, $n=1,2,\dots$ & $\Sp(n)$ & $n(2n+1)$ \\
\hline $\D_n$, $n=3,4,\dots$ & $\SO(2n)$ & $n(2n-1)$
\\
\hline
$\G_2$ & $\Aut(\oo)$ & $14$ \\
\hline $\F_4$ & $\Isom(\oo\pp^2)$ & $52$
\\
\hline $\E_6$ & $\Isom((\cc \otimes \oo)\pp^2)$ &$78$ \\
\hline $\E_7$ & $\Isom((\hh \otimes \oo)\pp^2)$& $133$
\\
\hline $\E_8$ & $\Isom((\oo \otimes \oo)\pp^2)$& $248$\\
\hline
\multicolumn{3}{|l|}{{\footnotesize the index $n$ denotes the rank}}\\
\hline
\end{tabular}
\end{table}
From theorem \cite{koll}.2.2, p.~13, we cite tables
\ref{PQKtable05}, \ref{PQKtable06}, \ref{PQKtable07} of maximal
connected subgroups (up to conjugation) of the classical Lie groups.
(By $\operatorname{Irr}_\rr$, $\operatorname{Irr}_\cc$,
$\operatorname{Irr}_\hh$ real, complex and quaternionic irreducible
representations are denoted. The tensor product ``$\otimes$'' of
matrix Lie groups is induced by the Kronecker product of matrices.)
\begin{table}[h]
\centering \caption{Maximal connected subgroups of $\SO(n)$}
\label{PQKtable05}
\begin{tabular}{|c@{\hspace{8mm}}| @{\hspace{8mm}}l@{\hspace{8mm}}|}
\hline subgroup & for \\
\hline\hline  $\SO(k)\times \SO(n-k)$& $1\leq k\leq n-1$ \\
\hline $\SO(p)\otimes \SO(q)$ & $pq=n$, $3\leq p\leq q$\\
\hline $\U(k)$ & $2k=n$ \\
\hline
$\Sp(p)\otimes \Sp(q)$ & $4pq=n$\\
\hline$\varrho(H)$ & $H$ simple, $\varrho\in \operatorname{Irr}_\rr(H)$, $\deg \varrho=n$\\
\hline
\end{tabular}
\end{table}

\begin{table}[h]
\centering \caption{Maximal connected subgroups of $\SU(n)$}
\label{PQKtable06}
\begin{tabular}{|c@{\hspace{8mm}}| @{\hspace{8mm}}l@{\hspace{8mm}}|}
\hline subgroup & for \\
\hline\hline  $\SO(n)$&  \\
\hline $\Sp(m)$ & $2m=n$\\
\hline $\mathbf{S}(\U(k)\times \U(n-k))$ & $1\leq k\leq n-1$ \\
\hline
$\SU(p)\otimes \SU(q)$ & $pq=n$, $p\geq 3$, $q\geq 2$\\
\hline$\varrho(H)$ & $H$ simple, $\varrho\in \operatorname{Irr}_\cc(H)$, $\deg \varrho=n$\\
\hline
\end{tabular}
\end{table}

\begin{table}[h]
\centering \caption{Maximal connected subgroups of $\Sp(n)$}
\label{PQKtable07}
\begin{tabular}{|c@{\hspace{8mm}}| @{\hspace{8mm}}l@{\hspace{8mm}}|}
\hline subgroup & for \\
\hline\hline  $\Sp(k)\times \Sp(n-k)$& $1\leq k\leq n-1$ \\
\hline $\SO(p)\otimes \Sp(q)$ & $pq=n$, $p\geq 3$ , $q\geq 1$\\
\hline $\U(n)$ &  \\
\hline$\varrho(H)$ & $H$ simple, $\varrho\in \operatorname{Irr}_\hh(H)$, $\deg \varrho=2n$\\
\hline
\end{tabular}
\end{table}
From \cite{bor}, p.~219, we cite table \ref{PQKtable08} of maximal
rank maximal connected subgroups.
\begin{table}[h]
\centering \caption{Maximal rank maximal connected subgroups}
\label{PQKtable08}
\begin{tabular}{|c@{\hspace{2mm}}| @{\hspace{2mm}}c@{\hspace{2mm}}|}
\hline ambient group & subgroup \\
\hline\hline  $\SU(n)$& $\mathbf{S}(\U(i)\times \U(n-i-1))$  for $i\geq 1$ \\
\hline $\SO(2n+1)$ & $\SO(2n)$, $\SO(2i+1)\times \SO(2(n-i))$  for $1\leq i\leq n-1$\\
\hline $\Sp(n)$ & $\Sp(i)\times \Sp(n-i)$ for $i\geq 1$, $\U(n)$  \\
\hline $\SO(2n)$ & $\SO(2i)\times \SO(2(n-i))$ for $i\geq 1$, $\U(n)$\\
\hline $\G_2$ & $\SO(4)$, $\SU(3)$\\
\hline
\end{tabular}
\end{table}
From table \cite{koll}.2.1 we cite subgroups of maximal dimension in
table \ref{PQKtable09}.
\begin{table}[h]
\centering \caption{Subgroups of maximal dimension}
\label{PQKtable09}
\begin{tabular}{|c@{\hspace{8mm}}| @{\hspace{8mm}}c@{\hspace{8mm}}|}
\hline ambient group & subgroup \\
\hline\hline  $\SU(n)$, $n\neq 4$& $\mathbf{S}(\U(1)\times \U(n+1))$ \\
\hline $\SU(4)$& $\Sp(2)$  \\
\hline
$\SO(n)$ & $\SO(n-1)$\\
\hline $\Sp(n)$, $n\geq 2$ & $\Sp(n-1)\times \Sp(1)$  \\
\hline $\G_2$ & $\SU(3)$\\
\hline
\end{tabular}
\end{table}

\vspace{3mm}

\label{page05}Let $H\In \Isom(M)$ be either $\s^1$ or $\zz_2$.
Consider the isotropy representation at an $H$-fixed-point $x\in M$
composed with the canonical projection $\Sp(1)\to \SO(3)$:
\begin{align*}
\varphi: H\hto{}\Sp(n)\Sp(1)\to \SO(3)
\end{align*}
Theorems \cite{danc}.4.4, p.~602, and \cite{danc}.5.1, p.~606,
(together with \cite{danc}, p.~600) show that the type of the
fixed-point component $F$ of $H$ around $x$ depends on the image of
$\varphi$: If $\varphi(H)=1$, the component $F$ is quaternionic for
$H\in \{\s^1,\zz_2\}$. If $\varphi(H)\neq 1$, the component $F$ is
locally K\"ahlerian for $H=\zz_2$ and K\"ahlerian for $H=\s^1$. A
result by Gray (cf.~\cite{gra}) shows that a quaternionic
submanifold is totally geodesic. Formula \cite{bess}.14.42b, p.~406,
then yields that the quaternionic components are again Positive
Quaternion K\"ahler Manifolds.

It is easy to see that the dimension of $F$ is exactly $2n$ if
$H=\zz_2$ and given that $F$ is locally K\"ahlerian. The dimension
of $F$ is smaller than or equal to $2n$ if $H=\s^1$ and provided
that $F$ is K\"ahlerian.

\vspace{3mm}

Let us finally state some cohomological properties of Positive
Quaternion K\"ahler Manifolds.

\begin{theo}[Cohomological properties]\label{PQKtheo05}
A Positive Quaternion K\"ahler Manifold $M$ satisfies:
\begin{itemize}
\item
Odd-degree Betti numbers vanish, i.e.~$b_{2i+1}=0$ for $i\geq 0$.
\item The identity
\begin{align*}
\sum_{p=0}^{n-1}\big(6p(n-1-p)-(n-1)(n-3)\big)b_{2p}=\frac{1}{2}
n(n-1)b_{2n}
\end{align*}
holds and specialises to
\begin{align}
\label{PQKeqn02}-1+3b_2+3b_4-b_6&=2b_8\\
\label{PQKeqn03}-4+5b_2+8b_4+5b_6-4b_8&=5b_{10}
\end{align}
in dimensions $16$ and $20$ respectively.
\item
A Positive Quaternion K\"ahler Manifold $M^{4n}\not\cong
\Gr_2(\cc^{n+2})$ is rationally
$3$-connected.
\item
The real cohomology algebra possesses an analogue of the
\linebreak[4]\emph{Hard-Lefschetz property}, i.e.~with the four-form
$u\in H^4(M,\rr)$ from above the morphism
\begin{align*}
L^k: H^{n-k}(M,\rr) \to H^{n+k}(M,\rr) \qquad
L^k(\alpha)=u^k \wedge \alpha
\end{align*}
is an isomorphism. In particular, we obtain
\begin{align*}
b_{i-4} \leq b_{i}
\end{align*}
for (even) $i\leq 2n$. A generator in top cohomology $H^{4n}(M)$ is
given by $u^n$. This defines a canonical orientation.
\end{itemize}
\end{theo}
\begin{prf}
The first point is proven in theorem \cite{sal1}.6.6, p.~163, where
it is shown that the Hodge decomposition of the twistor
space is concentrated in terms
$H^{p,p}(Z)$. The second item is due to
\cite{sal3}.5.4, p.~403.

The next item basically follows from theorem \cite{sal2}.5.5, p.~103
where it is proven that $b_2=0$ for $M^n\not\cong \Gr_2(\cc^{n+2})$.

The Hard-Lefschetz property of $M$ follows from the Hard-Lefschetz
property of the twistor space. Indeed, a Positive Quaternion
K\"ahler Manifold has this property with respect to $u$.
\end{prf}
So for a Positive Quaternion K\"ahler Manifold $M$ it is equivalent
to demand that $M$ be rationally $3$-connected---i.e.~to have that
$\pi_1(M)\otimes \qq=\pi_2(M)\otimes
\qq=\pi_3(M)\otimes\qq=0$---and to
require that $M$ be $\pi_2$-finite---i.e.~to suppose that
$\pi_2(M)<\infty$.


\section{Preparations}\label{sec02}

This section is devoted to a computation of several indices
$i^{p,q}$ in terms of the characteristic numbers of the complexified
tangent bundle $T_\cc M$ for a Positive Quaternion K\"ahler Manifold
$M$ of dimension $20$. That is, we compute
\begin{align*}
i^{p,q}=\langle \hat A(M)\cdot \ch(R^{p,q}), M\rangle
\end{align*}
with $R^{p,q}=\bigwedge_0^p E\otimes S^q H$ as usual. The formulas
relating these indices to other invariants are given in theorem
\ref{PQKtheo10}. Combining these equations with our computations
yields the fundamental system of equations we shall mainly be
concerned with in the following. It is linear in the characteristic
numbers of $M$.

We shall compute these indices in terms of characteristic classes
$u=-c_2(H)$ respectively $c_2, c_4,\dots, c_{10}$ of the bundles $H$
and $E$. Using the formula $\ch(E)=\sum_{i=1}^{10} e^{x_i}$ (for the
formal roots $x_i$), the analogue for the bundle $H$ and the fact
that Chern classes may be described as the elementary symmetric
polynomials in the formal roots one obtains easily:
\begin{align*}
\ch(H)=&2+u+\frac{u^2}{12}+\frac{u^3}{360}+\frac{u^4}{20160}+\frac{u^5}{1814400}\\
\ch(E)=&10-c_2+\frac{1}{12}\big(
c_2^2-2c_4\big)+\frac{1}{360}\big(-c_2^3+3c_2c_4-3c_6\big)
\\&+\frac{1}{20160}\big(c_2^4-4c_2^2c_4+2c_4^2+4c_2c_6-4c_8\big)
\\&+\frac{1}{1814400}\big(-5c_{10}-c_2^5+5c_2^3c_4-5c_2c_4^2-5c_2^2c_6+5c_4c_6+5c_2c_8\big)
\end{align*}
Now use the formula $T_\cc M=E\otimes H$ to successively compute the
Chern classes of the complexified tangent bundle and the Pontryagin
classes $p_i$ of $M$. Filling in these Pontryagin classes into the
characteristic series of the $\hat A$-genus yields
\begin{align*}
\hat A(M)=&1+\frac{1}{12}\big(c_2-5u\big)+\frac{1}{720}(3c_2^2-c_4-28c_2u+65u^2\big)
\\&+\frac{1}{60480}(10 c_2^3-9c_2c_4+2c_6-136 c_2^2 u +55 c_4 u + 570 c_2 u^2
 - 820 u^3 \big)\\&+\frac{1}{3628800}\big(21 c_2^4-34 c_2^2 c_4 +5 c_4^2 +13 c_2c_6 -3 c_8
-384 c_2^3 u \\&+ 409 c_2 c_4 u
-113 c_6 u + 2274 c_2^2 u^2-1060 c_4 u^2-5736 c_2 u^3 +5760 u^4\big)
\\&+\frac{1}{479001600}\big(90c_2^5-219 c_2^3c_4+87c_2c_4^2+109c_2^2c_6-32c_4c_6-43c_2c_8\\&+10 c_{10}
-2136 c_2^4 u+3990c_2^2c_4 u-675 c_4^2 u-1834 c_2c_6 u +525c_8 u
\\&+16524 c_2^3 u^2 -19740 c_2c_4 u^2+6155 c_6 u^2-57576 c_2^2
u^3+29935 c_4 u^3 \\&+98815 c_2 u^4 - 73985 u^5\big)
\end{align*}
We now compute the Chern characters of the exterior powers of $E$.
For this we use that the roots of $\bigwedge^k E$ are given by
$y_{i_1,\dots,i_k}=x_{i_1}+\dots + x_{i_k}$ for $1\leq i_1 < \dots <
i_k\leq 10$. It remains to compute the Chern characters of the
symmetric bundles, which can be done in a similar fashion e.g.~using
the roots given by the $y_{i_1,\dots,i_k}=x_{i_1}+ \dots+ x_{i_k}$
for $1\leq i_1\leq \dots \leq i_k\leq 10$.

This enables us to compute the following indices in terms of
characteristic numbers via \eqref{PQKeqn06}. (The ``indices"
$i^{p,q}$ with $p+q+5$ odd are to be regarded as formal expressions,
as they do not necessarily correspond to twisted Dirac operators.)
\begin{align*}
i^{0,0}=&\frac{1}{479001600}\big( 10c_{10}+90 c_2^5-32 c_4 c_6 -2136
c_2^4 u -675 c_4^2 u +525 c_8 u\\& +6155 c_6 u^2 +29935 c_4 u^3
-73985 u^4- 3 c_2^3(73c_4 -5508 u^2)+ c_2^2(109 c_6 \\&+3990 c_4 u -
57576 u^3) +c_2(87 c_4^2-43 c_8-1834 c_6 u -19740 c_4 u^2\\& +98815 u^4) \big)\\
i^{0,1}=& \frac{1}{239500800} \big(10 c_{10} + 90 c_2^5 -32 c_4 c_6
-750 c_2^4 u -345 c_4^2 u +327 c_8 u\\& -643 c_6 u^2 -22799 c_4
u^3+90817 u^5- 3c_2^3(73 c_4+1840 u^2) +c_2^2(109c_6\\&+1746 c_4 u
+50400 u^3) + c_2(87 c_4^2-43 c_8 -976 c_6 u +4284 c_4 u^2\\&
-116543
u^4\big)\\
 i^{0,3}=&\frac{1}{119750400}\big( 10 c_{10} + 90 c_2^5-32
c_4c_6+4794 c_2^4 u  + 975 c_4^2 u - 465 c_8 u\\& - 4075 c_6 u^2
+87025 c_4 u^3 + 310465 u^5 - 3 c_2^3(73 c_4 - 8368 u^2) \\&+
c_2^2(109 c_6 -7230 c_4 u -135456 u^3)
\\&+ c_2(87 c_4^2-43 c_8+ 2456 c_6 u-6540 c_4 u^2 -277055 u^4)\big)\\
i^{0,5}=&\frac{1}{79833600}\big(10 c_{10}+90 c_2^5-32 c_4 c_6 +14034 c_2^4 u +3175 c_4^2 u -1785 c_8 u \\&+74685 c_6 u^2  -1546255 c_4 u^3 +44944065 u^5+c_2^3(-219 c_4+498544 u^2)\\&+c_2^2(109 c_6-22190 c_4 u +6228704 u^3)\\&+c_2( 87 c_4^2-43 c_8+8176 c_6 u - 404740 c_4 u^2 +30037185 u^4)\big)\\
i^{0,7}=&\frac{1}{59875200}\big(10 c_{10} + 90 c_2^5-32 c_4 c_6+26970 c_2^4 u +6255 c_4^2 u - 3633 c_8 u\\& + 362357 c_6 u^2 -18291599 c_4 u^3 + 3830160577 u^5 -3 c_2^3(73 c_4\\& -682800 u^2)+c_2^2(109 c_6-43134 c_4 u + 61087200 u^3)\\&+c_2(87 c_4^2-43 c_8 +16184 c_6 u-1760556 c_4 u^2 +796656577 u^4)\big)\\
i^{1,0}=&\frac{1}{239500800}\big(-610 c_{10} +450 c_2^5 +1952 c_4 c_6 -9294 c_2^4 u - 49575 c_4^2 u\\& + 22425 c_8 u-149405 c_6 u^2 +690875 c_4 u^3 -369925 u^5 + c_2^3(-2481 c_4\\&+60576 u^2)+c_2^2(-4669 c_6 +43050 c_4 u -179904 u^3)\\&+c_2(4593 c_4^2-3317 c_8 +56104 c_6 u-224760 c_4 u^2 +113915 u^4))\\
\end{align*}
\begin{align*}
i^{1,2}=&\frac{1}{79833600}\big(-610 c_{10}+450 c_2^5+1952 c_4 c_6 +9186 c_2^4 u + 60425 c_4^2 u\\&- 43575 c_8 u + 8115 c_6 u^2 +750715 c_4 u^3 -693765 u^5-c_2^3(2481 c_4\\&+48544 u^2)-c_2^2(4669 c_6+39670 c_4 u +144704 u^3)\\& +c_2(4593 c_4^2-3317 c_8-101416 c_6 u+203800 c_4 u^2 +43515 u^4)\big)\\
i^{1,4}=&\frac{1}{47900160}\big(-610 c_{10}+450 c_2^5+1952 c_4c_6 +46146 c_2^4 u + 280425 c_4^2 u\\&- 175575 c_8 u -5810093 c_6 u^2 -31189765 c_4 u^3 -6917125 u^5 \\&- 3 c_2^3(827 c_4-333472 u^2)\\&-c_2(-4593 c_4^2+3317 c_8 +416456 c_6 u +3272904 c_4 u^2 +8410117 u^4)\\&- c_2^2(4669 c_6+4770(43 c_4 u -1504 u^3))\big)\\
i^{2,1}=&\frac{1}{5443200}\big(15070 c_{10}+90 c_2^5 - 2864 c_4 c_6
-246 c_2^4u +855 c_4^2 u + 21567 c_8 u \\&-153103 c_6 u^2-79439 c_4
u^3 +90817 u^5 -3 c_2^3(241 c_4 +2112 u^2)\\&+c_2^2(13 c_6+1146 c_4
u +21984 u^3)\\&
+c_2(1599 c_4^2+8369 c_8 - 7840 c_6 u +32784 c_4 u^2 +56257 u^4)\big)\\
i^{2,3}=& \frac{1}{2721600}\big(15070 c_{10}+90 c_2^5 -2864 c_4 c_6
+5298 c_2^4 u + 59775 c_4^2 u\\&+530535 c_8 u
+ 1506665 c_6 u^2 +60625 c_4 u^3 +310465 u^5 \\&+c_2^3(-723 c_4 +53088 u^2)\\&+c_2^2(13 c_6-36630 c_4 u +109728 u^3)\\&+c_2(1599 c_4^2+8369 c_8 -5848 c_6 u-307800 c_4 u^2 +414145 u^4)\big)\\
i^{3,0}=&\frac{1}{21772800}\big(-876370 c_{10}+450 c_2^5-38176
c_4c_6-7278 c_2^4 u -73575 c_4^2 u \\&
+1714425 c_8 u + 571315 c_6 u^2-293125 c_4 u^3- 369925 u^5\\&+c_2^3(-4497 c_4+28512 u^2)\\&+c_2^2(9347 c_6 + 55050 c_4 u -22848 u^3)\\&+ c_2(10641 c_4^2+15931 c_8-121112 c_6 u-159720 c_4 u^2-439045 u^4)\big)\\
i^{3,2}=&\frac{1}{7257600} \big(-876370 c_{10}+450 c_2^5-38176 c_4c_6+11202 c_2^4 u +190025 c_4^2 u \\&- 3765975 c_8 u -158205 c_6 u^2 -636485 c_4 u^3 -693765 u^5\\&-c_2^3(4497 c_4 +3808 u^2)\\&+c_2^2(9347 c_6-104470 c_4 u-225728 u^3)\\&+c_2(10641 c_4^2+15931 c_8+201368 c_6 u-126680 c_4 u^2-1062405 u^4)\big)\\
\end{align*}
\begin{align*}
i^{4,1}=&\frac{1}{7257600}\big(3509330 c_{10}+450 c_2^5-61216 c_4
c_6 + 282 c_2^4 u + 46275 c_4^2 u\\& +10275 c_8 u + 994705 c_6 u^2
+1245365 c_4 u^3+454085 u^5\\& - 3 c_2^3(1709 c_4+11376 u^2)\\&
+c_2^2(18977 c_6-15270 c_4 u +24672 u^3)\\& +c_2(12531 c_4^2-73079
c_8 +101488 c_6 u +228300 c_4 u^2 +799685 u^4)\big)\\
i^{5,0}=&\frac{1}{3628800}\big(-1415810 c_{10}+90 c_2^5-15488 c_4
c_6 -1254 c_2^4 u - 13275 c_4^2 u\\& - 1020075 c_8 u - 599905 c_6
u^2-314465 c_4 u^3 - 73985 u^5\\&-3 c_2^3(367 c_4-832
u^2)\\&+c_2^2(5191 c_6+10290 c_4 u + 11136 u^3)\\&-c_2(-2733
c_4^2+33097 c_8 +25816 c_6 u+3360 c_4 u^2+143105 u^4)\big)
\end{align*}
Using this information one may form the described linear system of
equations.

\vspace{3mm}

Now compute the \emph{Hilbert Polynomial} $f$ of $M$ in the
parameters $d$, $v$ and $i^{0,0}\in\qq$, i.e.~in the dimension of
the isometry group, the quaternionic volume and the $\hat A$-genus. The
Hilbert Polynomial $f$ on $M$ is given by
\begin{align*}
f(q)=\ind \Dc(S^q H)=\langle \hat A\cdot \ch(S^q
H),[M]\rangle=i^{0,q}
\end{align*}
and has degree $11$. We use the formula
\begin{equation*}
(H-2)^{\otimes m}=\sum_{j=0}^m (-1)^j \bigg( \tbinom{2m}{j} -
\tbinom{2m}{j-2}\bigg) S^{m-j}H
\end{equation*}
resulting from the Glebsch-Gordan formula. The leading term of the
power series of the $\hat A$-genus is $1$ and the first non-zero
coefficient of the power series $\ch(H-2)^m$  lies in degree $m$.
Thus we obtain that $\langle\hat A(M) \cdot \ch(H-2)^{\otimes
5},[M]\rangle=u^5$ and all the higher terms vanish. Combining this
with theorem \ref{PQKtheo10}, i.e.~with $f(0)=i^{0,0}$,
$f(1)=f(3)=0$, $f(5)=1$ and $f(7)=d$ permits us to compute the
following identities.
\begin{align*}
f(0)&=i^{0,0} \\ f(1)&=0 \\ f(2)&= \frac{-2816 + 128 d -360448 i^{0,0}- 7v}{229376}\\
f(3)&=0 \\ f(4)&=\frac{269568 - 7040 d +4685824 i^{0,0}+ 273 v}{1146880} \\ f(5)&=1 \\
f(6)&=\frac{228096 + 18304 d -2342912 i^{0,0}- 273 v}{114688} \\
f(7)&=d \\
f(8)&=\frac {13(-143616+35200 d+ 3063808 i^{0,0}+595 v)}{114688}\\
f(9)&=\frac{1}{140}(-10692 + 1760 d+262144 i^{0,0} +63 v)\\
f(10)&=\frac{13(-4333824 + 598400 d +116424704 i^{0,0}+ 33915 v)}{229376} \\
f(11)&=\frac{1}{14}(-9152+1144 d+262144 i^{0,0}+ 91 v)
\end{align*}
The Hilbert polynomial has degree $11$ and thus can be computed from
these values.

From theorem \cite{sem1}.1.1, p.~2, we are given the formula
\begin{equation*}
0\leq f_M(5+2q)\leq f_{\hh\pp^5}(5+2q)={11+2q\choose 11}
\end{equation*}
for $q \in \nn_0$. So we may compute for each $q$ a lower and an
upper bound for $i^{0,0}$---depending on $d$ and $v$. Unfortunately,
with $q$ growing, these bounds seem to become worse so that we use
low values of $q$---i.e.~$q=3$ respectively $q=2$---to obtain:
\begin{align}\label{PREPeqn01}
\frac{1}{14}(-9152+262144 i^{0,0} + 1144 d+91 v)&\leq 12376 \\
\label{PREPeqn02} \frac{1}{140}(-10692+262144 i^{0,0}+ 1760 d+63
v)&\geq 0
\end{align}

\vspace{3mm}

Let us now compute further relations involving the $\hat A$-genus of
a Positive Quaternion K\"ahler Manifold $M$. We adapt lemma
\cite{sal1}.7.6, p.~169, to dimension $20$ and plug in the
expression $-(2c_2-10u)$ for the first Pontryagin class $p_1$:
\begin{align}\label{PREPeqn03}
8 u^5-p_1 u^4 \geq 0 \Leftrightarrow 8u^5+(2 c_2-10 u)u^4\geq 0
\Leftrightarrow c_2u^4-u^5\geq 0
\end{align}
From the solution of the fundamental system of equations we cite
\begin{align*}
c_2u^4=-\frac{81}{70}+\frac{3 d}{28}+\frac{1536 \hat
A(M)[M]}{35}-\frac{31 u^5}{5}
\end{align*}
Combining this with formula \eqref{PREPeqn03} yields
\begin{align}\label{PREPeqn05}
-\frac{81}{70}+\frac{3d}{28}-\frac{36 u^5}{5}+\frac{1536 \hat
A(M)[M]}{35} \geq 0
\end{align}
Recall that $1\leq v \leq 1024$. So for $d=0$ the equation becomes
\begin{align}\label{PREPeqn04}
-\frac{81}{70}-\frac{36}{5\cdot 1024} +\frac{1536 \hat A(M)[M]}{35}
\geq 0
\end{align}


\section{Special cases and the proof of theorem B}\label{sec03}

This section will combine further observations with the proof of
theorem \ref{theoB}. Indeed, we shall deal with each
dimension---i.e.~$\dim M\in \{20,24\}$---in theorem \ref{theoB}
separately thereby proving slightly more general assertions. The
theorem itself is then a combination of corollary \ref{RECcor01} and
theorem \ref{RECtheo02}. As we already remarked, in dimension $28$
there is an exceptional Wolf space, which makes further
generalisation more difficult. Nonetheless, as we were told by
Gregor Weingart, a similar recognition theorem---which also
identifies the exceptional Wolf space $\F_4/{\Sp(3)\Sp(1)}$---seems
to be possible.

Foremost we recall
\begin{theo}\label{PQKtheo09}
Let $M$ be a Positive Quaternion K\"ahler Manifold. If $8\neq \dim M\leq
16$ and $b_4(M)=1$, then $M$ is homothetic to the quaternionic
projective space.
\end{theo}
\begin{prf}
See theorem \cite{sal2}.2.1.ii, p.~89.
\end{prf}
Before generalising this theorem, we shall reconsider the problem in
dimension $16$ and we shall illustrate the used methods by pointing
out certain additional results.

\subsection{Dimension $16$}

In dimension $16$ the relation on Betti numbers given in
\eqref{PQKeqn02} of theorem \ref{PQKtheo05} together with the
Hard-Lefschetz property (cf.~theorem \ref{PQKtheo05}) have the
following consequence: If $b_4=1$ and if we assume $M$ to be
rationally $3$-connected (cf.~\ref{PQKtheo05}), we obtain
$b_0=b_4=b_8=b_{12}=b_{16}=1$ with all the other Betti numbers
vanishing. So necessarily every Pontryagin class is a multiple of
the corresponding power of the form $u$. This motivates the
following slight improvement.
\begin{prop}\label{RECprop01}
If each of the Chern classes $c_i$ of the bundle $E$ over a
$16$-dimensional Positive Quaternion K\"ahler Manifold is a (scalar)
multiple of the corresponding power of $u$, then the manifold
already is homothetic to $\hh\pp^4$.
\end{prop}
\begin{prf}
We form a linear system of equations as we did in section \ref{sec02}.
By assumption we may now replace every Chern class $c_i$ by some
$x_iu^{n_i}$ for $x_i \in \rr$.

If one focuses on the case $b_2=0$ (cf.~\ref{PQKtheo05}), the system
of equations can be solved and it yields $d=55$, $b_4=b_8=1$,
$b_6=0$, $u^4=1$ (with all the factors $x_i$ equal to one). We then
observe that in dimension $55$ only semi-simple Lie groups of rank
at least $5$ appear. Theorem \ref{PQKtheo07} then yields the
assertion; i.e.~the isometry group becomes very large and permits to
identify $M$ as the quaternionic projective space.

If one does not assume $b_2=0$, the list of possible configurations
for $(d,b_2,b_4,b_6,b_8,u^4)$ becomes a little larger. However, the
configuration from above remains the only one with integral $d\in
\zz$.
\end{prf}
Assume $b_2=0$. Then the same proof works if one only requires $c_2$
and $c_4$ to be scalar multiples of $u$ respectively $u^2$. In this
case a numerical solving procedure leads to six different solutions
of which the only one with an integral value for $d$ is the
requested one---as in proposition \ref{RECprop01}.

Focussing on the case that only $c_2$ is a multiple $x\in \rr$ of
$u$ leads to the two equations
\begin{align}\label{RECeqn11}
d&=7+\frac{v}{6}+\frac{vx}{48}\\
\label{RECeqn12}
b_4&=\frac{783}{2}-\frac{7}{8}v-\frac{9}{16}vx-\frac{11}{128}vx^2
-\frac{1}{512}vx^3
\end{align}
where $v=(4u)^4$ is the quaternionic volume. The element $x$ is
integral by the same reasoning as in the original proof, i.e.~the
proof of theorem \cite{sas}.5.1, p.~62. In \cite{brun} it is proven
that $i^{1,n+1}\leq 0$. In the survey article \cite{sal2}, p.~117,
it is suggested that this index vanishes unless $M$ is the
quaternionic projective space. In the following we assume the
vanishing of $i^{1,5}$ in the case $M\neq \hh\pp^4$, which produces
\begin{align}
\label{RECeqn13}
d=&\frac{7(304+56x+3x^2)}{16+20x+3x^2}\\
\label{RECeqn14}
b_4=&-\frac{27(1280-304x-40x^2+7x^3)}{8(16+20x+3x^2)}\\
\label{RECeqn15}
b_6=&\frac{1}{36}\big(3289-294x+63x^2-6(c_4u^2)x^2-\frac{53200}{16+20x+3x^2}
\\&\nonumber-\frac{93548x}{16+20x+3x^2}\big)\\
\label{RECeqn16} b_8=&\frac{1}{144}\bigg(
14410-1113x-(126x^2-12(c_4u^2))x^2-
\frac{1163680}{16+20x+3x^2}\\&\nonumber+\frac{14728x}{16+20x+3x^2}\bigg)\\
\label{RECeqn17}
c_4^2=&\frac{1}{16}\bigg(-3546+567x-378x^2+44(c_4u^2)x^2
-\frac{821520}{16+20x+3x^2}\\&\nonumber+\frac{445032x}{16+20x+3x^2}\bigg)
\end{align}
The only integral solution for $8\leq d<55=\dim \Sp(5)$ (cf.~theorem
\ref{PQKtheo07}) is given by $x=4$ and $d=28$. Then we directly
obtain $b_4=3$ by \eqref{RECeqn14} and also $v=84$ by
\eqref{RECeqn11}. Indeed, by the relations on Betti numbers in
\ref{PQKtheo05} only two possibilities for $(b_4,b_6,b_8)$ remain,
namely $(3,0,4)$ or $(3,2,3)$. Equations \eqref{RECeqn15} and
\eqref{RECeqn16} yield $c_4u^2=\frac{27}{32}$ in the first case. By
theorem \ref{PQKtheo08} we may use the positive definiteness of the
generalised intersection form $Q$ to see $0\leq Q(c_4,c_4)=c_4^2$.
Yet, in the case $(b_4,b_6,b_8)=(3,2,3)$ we obtain the contradiction
$c_4^2=-\frac{75}{16}$ by \eqref{RECeqn17}. As a consequence, we
have the following theorem:
\begin{theo}\label{theo07}
If $M$ is a rationally $3$-connected $16$-dimensional Positive
Quaternion K\"ahler Manifold with $i^{1,5}=0$ and if the class $c_2$
is a scalar multiple of $u$, then either $M\cong \hh\pp^4$ or the
datum $(d,v,b_4,b_6,b_8)=(28,84,3,0,4)$ is exactly the one of
$\widetilde \Gr_4(\rr^8)$.
\end{theo}
\vproof We remark that the property that $c_2$ is a multiple of $u$
seems to be a special feature of $\widetilde\Gr_4(\rr^{n+4})$ for
$n=4$ among the infinite series of Wolf spaces other than the
quaternionic projective space.

\subsection{Dimension $20$}

The following consequence is as simple as it is astonishing.
\begin{lemma}\label{PQKlemma01}
Let $M$ be rationally $3$-connected of dimension $20$ with $b_4\leq
5$. Then it holds:
\begin{align}\label{eqn14}
b_6=b_{10} \quad \vee \quad
(b_4,b_8)\in\{(1,1),(2,3),(3,5),(4,7),(5,9)\}
\end{align}
\end{lemma}
\begin{prf}
By assumption $b_2=0$. Equation \eqref{PQKeqn03} becomes
\begin{align*}
4(2b_4-b_8-1)=5(b_{10}-b_6)
\end{align*}
where the right hand side is non-negative due to Hard-Lefschetz
(cf.~\ref{PQKtheo05}). Hence the term $2b_4-b_8-1$ must either be a
positive multiple of $5$ or zero. Since also $b_8\geq b_4$, the
first case may not occur for $b_4\leq 5$. Thus it holds that
$b_6=b_{10}$. The existence of the form $0\neq u\in H^4(M)$ shows
that $b_4\geq 1$.
\end{prf}

Let us now prove the analogue of proposition \ref{RECprop01} in
dimension $20$. Again, without restriction, we focus on rationally
$3$-connected Positive Quaternion K\"ahler Manifolds $M^{20}$.
\begin{theo}\label{RECtheo01}
If each of the Chern classes $c_{2i}$ of $E$ over $M^{20}$ is a (scalar)
multiple of the corresponding power of $u$, then $M^{20}\cong \hh\pp^5$.
\end{theo}
\begin{prf}
We proceed as before in dimension $16$; i.e.~we solve the system of
equations (cf.~section \ref{sec02}) setting all characteristic
classes to rational multiples of a suitable power of the form $u$.
As a result we obtain $b_4=b_8=1$, $b_6=b_{10}=0$, $u^5=1$, all the
scalars are $1$ and $d=78$. However, such a  large isometry group
can only occur for the quaternionic projective space---cf.~theorem
\ref{PQKtheo07}.
\end{prf}
We remark that for this proof to work we do not need that $c_6$ is a
scalar multiple of $u^3$. We shall now use the observation we stated
in lemma \ref{PQKlemma01} to finish the reasoning.
\begin{cor}\label{RECcor01}
If $b_4(M^{20})=1$, then $M^{20}\cong \hh\pp^5$.
\end{cor}
\begin{prf}
Lemma \ref{PQKlemma01} tells us that $b_4=1$ implies $b_8=1$. By
Poincar\'e Duality we may conclude that all the $c_{2i}\in
H^{4i}(M)$ satisfy the condition from theorem \ref{RECtheo01}.
\end{prf}
Observe that clearly by this corollary we have ruled out an infinite
number of possible configurations of Betti numbers, since it
automatically follows that $b_4=1$ not only implies $b_8=1$ but also
$b_6=b_{10}=0$.

Observe that one may prove this corollary in a slightly different
way: Assume only that $c_2$ and $u$ are scalar multiples and do the
same for monomials containing $c_2$ and $u$. Use the equations with
the additional information $b_4=1$ (and $b_2=0$) and the result
follows directly.


\subsection{Dimension $24$}

We apply similar techniques as before to prove
\begin{theo}\label{RECtheo02}
A $24$-dimensional Positive Quaternion K\"ahler Manifold $M^{24}$ with
$b_4(M)=1$ is homothetic to $\hh\pp^6$.
\end{theo}
\begin{prf}
Again we replace $c_2$ by a scalar multiple of $u$ and do the same
for all the monomials containing $c_2$ as a factor. This and the
additional information $b_4=1$ (respectively $b_2=0$) simplifies the
system of equations (cf.~\ref{PQKtheo10}) that we build up as we did
for dimensions $16$ and $20$. Solving it yields a list of
possibilities, which we run through in order to isolate the one
leading to $\hh\pp^6$:

We may rule out the first solution, since it yields
$d=\frac{244724}{2891}\not\in \zz$ for the dimension of the isometry
group. The second solution gives $d=105=\dim \Sp(7)$. Thus we
directly know that $M\cong \hh\pp^6$ in this case---cf.~theorem
\ref{PQKtheo07}. So what is left is to rule out the following
configuration of solutions, which is marked by
\begin{align*}
d=&\frac{7937019926774969}{402874803650560}x_1-\frac{19592196959405797}{2417248821903360}x_2^2+\frac{457452279096536909}{2417248821903360}
\\&-\frac{263256496233}{805749607301120}x_3^6+\frac{1176648936457}{402874803650560}x_4^5+\frac{14142811929437}{161149921460224}x_5^4\\&-\frac{282904843313851}{604312205475840}x_6^3
\end{align*}
where each $x_i$ is a root of
\begin{align*}
&29223x^7-358275x^6-6960405x^5+67759961x^4+579930789x^3\\&-4142432537x^2-9711667063x+33284884867
\end{align*}
Numerically, the roots of this polynomial are given by
\begin{align*}
&2.156753156, 7.720829360, 11.12408307, 15.23992325,\\&
-10.15093795+2.570319306i, -3.679678028,\\&
-10.15093795-2.570319306i
\end{align*}
A computer-based check on all the possible combinations now shows
that there are no integral solutions for $d$ in all these cases. So
we are done.
\end{prf}


\section{Properties of interest}\label{sec04}

In this section we shall show that under slight assumptions some
surprising results on the degree of symmetry of $20$-dimensional
Positive Quaternion K\"ahler Manifolds $M$ are obtained. This will
lead us to the existence of large isometry groups under mild
assumptions thereby making a first step towards the proof of theorem
\ref{theoA}.

\begin{prop}\label{INTprop01}
Unless $M$ admits an isometric $\s^1$-action, the $\hat
A$-genus of $M$ is restricted by
\begin{align*}
0.0321350097< \hat A (M)[M]< 0.6955146790
\end{align*}
\end{prop}
\begin{prf}
The upper bound clearly results from equation \eqref{PREPeqn01} when
substituting in the extremal value $v=1$. For the lower bound we
form the linear combination
\begin{align*}
\frac{1}{448}(-1053+136 d+32768 i^{0,0})\geq 0
\end{align*}
out of equations \eqref{PREPeqn02} and \eqref{PREPeqn05}. The result
follows from setting $d=0$.
\end{prf}

Let us now use the fact that the terms $f(5+2q)=i^{0,5+2q}$ (for
$q\geq 0$) are indices of the twisted Dirac operator
$\Dc(S^{5+2q}H)$; i.e.~in particular they are integral. This leads
to congruence relations for the dimension of the isometry group and
the quaternionic volume.
\begin{theo}\label{INTtheo01}
A $20$-dimensional rationally $3$-connected Positive Quaternion
K\"ahler Manifold with $\hat A(M)[M]=0$ satisfies
\begin{align*}
d\equiv 1 \mod 7 \qquad\qquad\textrm{and}\qquad\qquad v\equiv 4 \mod
20
\end{align*}
\end{theo}
\begin{prf}
We use the Hilbert Polynomial $f$ of $M$. Thus, under the assumption
that $\hat A(M)[M]=0$ we obtain
\begin{equation*}
\zz \ni i^{0,9}=f(9)=\frac{1}{140} (-10692 + 1760 d + 63 v).
\end{equation*}
This implies that
\begin{align*}
&-10692 + 1760 d +63 v \equiv 0 \mod 140 \\
\iff &88 + 80 d +63 v \equiv 0\mod 140 \\
\iff &(d \equiv 1 \mod 7) \vee (v\equiv 4 \mod 20)
\end{align*}
\end{prf}
\begin{rem}\label{INTrem01}
\begin{itemize}
\item
Any computation of further indices seems to result in the fact that
only denominators appear that divide $2^2 \cdot 5 \cdot 7$. Since
$v_{\hh\pp^5}=1024$ and since $v_{\widetilde\Gr_4(\rr^9)}=264$, we
see that the relations found in the theorem are the only ones that
may hold on the quaternionic volume when focussing on congruence
modulo $m$ for $m|140$.
\item
The dimension $d_{\hh\pp^n}$ of $\hh\pp^n$ is given by $(n+1)(2n+3)$
with
\begin{align*}
d_{\hh\pp^n}\equiv (n+1)((n+1)+(n+2))\equiv 1 \mod
n+2.
\end{align*}
The dimension $d_{\widetilde \Gr_4(\rr^{n+4})}$ of $\widetilde
\Gr_4(\rr^{n+4})$ is given by
\begin{align*}
d_{\widetilde \Gr_4(\rr^{n+4})}=
\begin{cases}
\frac{1}{2}(n+3)(n+4)&\textrm{for $n$ odd}\\
\frac{1}{2}(n+3)(n+4)&\textrm{for $n$ even}
\end{cases}
\end{align*}
So in any case we have $d_{\widetilde \Gr_4(\rr^{n+4})}\equiv 1 \mod n+2$.
Thus in dimensions without exceptional Wolf spaces---for these it is
not true---by the main conjecture \ref{PQKconj01} it should hold on
a rationally $3$-connected Positive Quaternion K\"ahler Manifold
that the dimension of the isometry group is congruent $1$ modulo
$n+2$.
\end{itemize}
\end{rem}
\begin{cor}\label{INTcor01}
A $20$-dimensional Positive Quaternion K\"ahler Manifold satisfying
$\hat A(M)[M]=0$ and possessing an isometry group of dimension
greater than $36$ is isometric to the complex Grassmannian or the
quaternionic projective space.
\end{cor}
\begin{prf}
Again we assume the manifold to be rationally
$3$-connected. There are no compact Lie groups in dimensions $43$,
$50$, $57$, $64$ and $71$ with rank smaller than or equal to $5$.
Thus theorems \ref{PQKtheo07} and \ref{INTtheo01} yield the result.
\end{prf}

\vspace{3mm}

We shall now prove the existence of $\s^1$-actions on
$20$-dimensional Positive Quaternion K\"ahler Manifolds which
satisfy some slight assumptions. For this we use the definiteness of
the intersection form to establish
\begin{lemma}\label{INTlemma01}
We have:
\begin{align*}
u^5\geq 0 \qquad\qquad c_2^4u\geq 0 \qquad\qquad  c_2^2u^3\geq 0
\qquad\qquad c_4^2 u\geq 0
\end{align*}
More generally, the same holds for
\begin{align*}
(k c_2^2 + l c_2 u + m u^2 + n c_4)^2 u\geq 0
\end{align*}
with $k,l,m,n \in \rr$.
\end{lemma}
\begin{prf}
Recall the generalised intersection form $Q$ from theorem
\ref{PQKtheo08}. All the classes $y$ from the assertion may be
written as $y=Q(x,x)$ for some $x\in H^r(M)$ with $r\in \{4,8\}$.
However, the intersection form $Q$ is positive definite in degrees
divisible by $4$ and it results that $Q(x,x)\geq 0$.
\end{prf}

Lemma \ref{INTlemma01} yields that $c_2^2u^3\geq 0$, which
translates to
\begin{align*}
\frac{495392 - 14240 d - 35651584 i^{0,0} - 1120 i^{1,6} + 707
v}{35840}\geq 0
\end{align*}
after solving the linear system of equations on indices (cf.~section
\ref{sec02}) and substituting in the special solution for $i^{1,6}$.
So we obtain:
\begin{align*}
i^{1,6} \leq \frac{ 495392 - 14240 d - 35651584 i^{0,0} + 707
v}{1120}
\end{align*}
which already shows that for extremal values of $v$ and $d$,
i.e.~$v=1024$ and $d=0$, the index $i^{1,6}$ becomes very small.

Now consider the term $(c_2u+m u^2)^2u$ together with the solution
from the system of equations (cf.~section \ref{sec02}) and the
solution for $i^{1,6}$. By lemma \ref{INTlemma01} we have $(c_2u+m
u^2)^2u\geq 0$. Suppose $d=0$. This has the consequence that
\begin{align*}
i^{1,6}\leq & \frac{1}{1120}(495392-35651584 i^{0,0} - 82944 m + 3145728
i^{0,0} m +707 v\\& - 434 m v +35 m^2 v)
\end{align*}
for all $m\in \rr$. The right hand side is a parabola in $m$.
Determine the apex of this parabola as $m_0=\frac{41472 - 1572864
i^{0,0} + 217 v}{35 v}$, put it into the inequality and obtain:
\begin{align*}
i^{1,6}\leq -\frac{1}{4900 v}&(214990848+309237645312
(i^{0,0})^2+82516 v+2793 v^2\\&+131072 i^{0,0}(-124416+539 v))
\end{align*}
The right hand side is a function in $i^{0,0}$ and $v$ which has no
critical point in the interior of the square $\big[0.0321350097,
0.695514790\big]\times \big[1, 1024\big]$---cf.~proposition
\ref{INTprop01}. Thus its maximum lies in the boundary of the
square. A direct check reveals that on the border of the square the
function is decreasing monotonously in $i^{0,0}$ for
$v\in\{1,1024\}$. Analogously, we see that for
$i^{0,0}=0.0321350097$ the function has the only maximum $-549.348$
for $v=61$ and for $i^{0,0}=0.695514790$ it is increasing in
$[1,1024]$. So it takes its maximum $-549.348$ on the square in
$v=61$ and $i^{0,0}=0.695514790$. So, in particular, we obtain the
following theorem:
\begin{theo}\label{INTtheo03}
A $20$-dimensional Positive Quaternion K\"ahler Manifold with
\begin{align*}
i^{1,6}\geq -549
\end{align*}
admits an effective isometric $\s^1$-action.
\end{theo}
\vproof As we remarked already for $M\not \cong
\hh\pp^n$, the index $i^{1,6}$ is smaller or equal to zero and it is
conjectured to equal zero. On $\hh\pp^n$ it equals
$i^{1,n+1}=n(2n+3)$.

\vspace{3mm}

Next we shall link the existence of an isometric $\s^1$-action to
the Euler characteristic.
\begin{theo}\label{INTtheo02}
A $20$-dimensional Positive Quaternion K\"ahler Manifold $M$ with
Euler characteristic restricted by
\begin{align*}
\chi(M)&<16236 \intertext{admits an effective isometric
$\s^1$-action. The same holds if the Betti numbers of $M$ satisfy
either} b_4-\frac{b_6}{4}&< 842.5 \intertext{or} \frac{59
b_4}{3}-\frac{25 b_6}{4}&<  3027.93 \intertext{or---as a combination
of both inequalities---if} b_4&\leq 3381
\end{align*}
\end{theo}
\begin{prf}
The theorem is trivial for $\Gr_2(\cc^{7})$. Thus we assume $M$ to
be rationally $3$-connected. We give a proof by contradiction and
assume $d=\dim \Isom(M)=0$. We shall choose special values for
$k,l,m,n$ from lemma \ref{INTlemma01}. These coefficients determine
an element $y$. We shall obtain the contradiction $Q(y,y)<0$ under
the assumptions from the assertion.

Use the linear combination in corollary \ref{INTlemma01} with
coefficients $n = -0.168$, $m=4.99$, $k=-n$, $l=-2\sqrt{-mn-18 n^2}$
under the assumption of $d=0$, $b_2=0$ together with our solution to
the system of indices (cf.~section \ref{sec02}) and the relation
on Betti numbers \eqref{PQKeqn03}. This results in the formula
\begin{align*}
19.9668 + 0.254016 b_4 - 0.063504 b_6 - 9835.62 i^{0,0} +
  0.0801763 v\geq 0
\end{align*}
(where coefficients are rounded off.) Substituting in the lower
bound $i^{0,0}\geq 0.0321350097$ from proposition \ref{INTprop01}
and the upper bound $v=1024$ yields $b_4-\frac{b_6}{4}\geq 842.468$.
This contradicts our assumption $b_4-\frac{b_6}{4} < 842.5$. Thus we
obtain $d\neq 0$.

\vspace{3mm}

The second formula involving Betti numbers results from similar
arguments with coefficients $l=22k$, $n=-\frac{3 l}{22}$,
$m=-\frac{239 n}{3}$, $n=1$. This yields
\begin{align*}
-467.202 + \frac{59 b_4}{3}-\frac{25 b_6}{4} - 104.025 i^{0,0} +
4.72222 i^{1,6} - 0.439931 v \geq 0
\end{align*}
Once more we assume there is no $\s^1$-action on $M$. Thus theorem
\ref{INTtheo03} gives us $i^{1,6}\leq -551$. So in this case we
additionally substitute the other known bounds $i^{0,0}\geq
0.0321350097$ and $v\geq 1$. This eventually yields that $\frac{59
b_4}{3}-\frac{25 b_6}{4}\geq 3027.93$ contradicting our assumption.

\vspace{3mm}

Assume $d=0$. Thus from the previous two relations on Betti numbers
we compute
\begin{align*}
25\cdot b_4-\frac{59 b_4}{3}&\geq 25\cdot842.5-  3027.93
\Leftrightarrow b_4\geq 3381.48
\end{align*}
Hence, whenever $b_4\leq 3381$, we obtain a contradiction and $d\neq
0$.

\vspace{3mm}

The result on the Euler characteristic results from the formula
$b_4\leq 3381$ and the Hard-Lefschetz property by a computer-based
check on all possible configurations of $(b_4,b_6,b_8,b_{10})$---in
a suitable range---that satisfy relation \eqref{PQKeqn03}. That is,
we start with $b_4=3382$ and figure out the configuration of
$(b_4,b_6,b_8,b_{10})$---satisfying all the properties from theorem
\ref{PQKtheo05}---with smallest Euler characteristic. This
configuration is given by
\begin{align*}
(b_4,b_6,b_8,b_{10})=(3382,0,3383,2704)
\end{align*}
and Euler characteristic $\chi(M)=16236$. So whenever
$\chi(M)<16236$ we necessarily have $b_4\leq 3381$. The result
follows by our previous reasoning.
\end{prf}
\begin{theo}\label{INTtheo04}
Let $M^{20}\not\in \{\hh\pp^5,\Gr_2(\cc^7)\}$ be a (rationally
$3$-connected) Positive Quaternion K\"ahler Manifold with $\hat
A(M)[M]=0$. Then it holds:
\begin{itemize}
\item
The dimension $d$ of the isometry group of $M$ satisfies
\begin{align*}
d\in \{15,22,29,36\}
\end{align*}
\item
The pair $(d,v)$ of the dimension of the isometry group and the
quaternionic volume is one of
\begin{align*}
&(15,4), (15,24), (15,44), (15, 64), (22,24), \dots , (22, 164),
\\&(29,24), \dots, (29, 264), (36, 24), \dots (36,384)
\end{align*}
where $v$ increases by steps of $20$.
\item
The connected component $\Isom_0(M)$ of the isometry group of $M$ is
as given in table \ref{INTtable01} up to finite coverings.
\begin{table}[h]
\centering\caption{Possible isometry groups \label{INTtable01}}
\begin{tabular}{|@{\hspace{4mm}}c@{\hspace{4mm}} |@{\hspace{4mm}}c@{\hspace{4mm}}|}
\hline $\dim \Isom(M)$ & type of $\Isom_0(M)$ up to finite coverings
\\
\hline\hline&\\[0mm] $15$ & $\SO(6)$, $\G_2\times \s^1$, $\SO(4)\times \SO(4)\times
\SO(3)$,\\& $\Sp(2)\times \Sp(1)\times \s^1\times \s^1$,
$\SU(3)\times \SO(4)\times \s^1$\\[4mm]
$22$ & $\Sp(3)\times \s^1$, $\SO(7)\times \s^1$, $\G_2\times \SU(3)$\\[4mm]
$29$ & $\SO(8)\times \s^1$, $\SO(6)\times \G_2$, $\G_2\times \G_2\times \s^1$,\\& $\SO(7)\times \SU(3)$, $\Sp(3)\times \SU(3)$ \\[4mm]
$36$ & $\SO(9)$, $\Sp(4)$\\[4mm]
 \hline
\end{tabular}
\end{table}
\end{itemize}
\end{theo}
\begin{prf}
Recall equation \eqref{PREPeqn05}
\begin{align}
\nonumber -\frac{81}{70}+\frac{3d}{28}-\frac{36 u^5}{5}+\frac{1536
\hat A(M)[M]}{35} \geq 0 \intertext{and set the $\hat A$-genus to
zero. This results in} \label{INTeqn01}
-\frac{162}{35}+\frac{3d}{7}-\frac{144 v}{1024\cdot 5}\geq 0
\end{align}
From computations with the Hilbert Polynomial we know that $d\equiv
1 \mod 7$ and $v\equiv 4 \mod 20$ for the dimension of the isometry
group and the quaternionic volume---cf.~theorem \ref{INTtheo01}.
Since $\frac{3d}{7}-\frac{162}{35}<0$ for all values of $d$ with $d
\equiv 1 \mod 7$ and $d<15$, we may thus rule them all out.

\vspace{3mm}

Now use the classification of compact Lie groups. The connected
component of the identity of $\Isom(M)$ permits a finite covering by
a product of a semi-simple Lie group and a torus (cf.~section
\ref{sec01}). Now we figure out all those products $G$ of simple Lie
groups and tori that satisfy
\begin{itemize}
\item $\dim G\equiv 1 \mod 7$,
\item
$15\leq \dim G\leq 36$ (cf.~corollary \ref{INTcor01}),
\item
$\rk G\leq 5$. (By theorem \ref{PQKtheo07} we know that $\rk G\geq
6$ already implies $M\cong\hh\pp^5$, as $M$ is rationally
$3$-connected.)
\end{itemize}
Congruence classes modulo $7$ of the dimensions of the relevant
different types of Lie groups are as described in table
\ref{INTtable02}.
\begin{table}[h]
\centering \caption{Dimensions modulo $7$} \label{INTtable02}
\begin{tabular}{|r@{\hspace{5mm}}| @{\hspace{5mm}}r@{\hspace{5mm}}|
@{\hspace{5mm}}r@{\hspace{5mm}}|@{\hspace{5mm}}
r@{\hspace{5mm}}|@{\hspace{5mm}}r@{\hspace{5mm}}|@{\hspace{5mm}}r@{\hspace{5mm}}|}
\hline type & $n=1$ & $n=2$ &$n=3$ & $n=4$ & $n=5$ \\
\hline\hline $\A_n$ & $3$ & $1$ & $1$ & $3$ & $0$\\
\hline
$\B_n$ & $3$ & $3$ & $0$ & $1$ & $6$\\
\hline
$\C_n$ & $3$ & $3$ & $0$ & $1$  & $6$\\
\hline
$\D_n$ & $1$  & $6$  & $1$ & $0$ & $3$\\
\hline
\multicolumn{6}{|c|}{{$\dim \G_2 \equiv 0 \mod 7$}}\\
\hline
\multicolumn{6}{|c|}{{$\dim \F_4 \equiv 3 \mod 7$}}\\
\hline
\end{tabular}
\end{table}
The list of Lie groups $G$ then is given as in the assertion.

\vspace{3mm}

Let us finally establish the list of pairs $(d,v)$. For this we
apply equation \eqref{INTeqn01} once more. Additionally, note that
the total deficiency
\begin{align*}
\Delta=2n+1+2v-d\geq 0
\end{align*}
is non-negative (cf.~\cite{herr2}, p.~208), which translates to
\begin{equation*}
11+2v-d\geq 0
\end{equation*}
for dimension $20$. These conditions already reduce the list of
pairs $(d,v)$ due to the following restrictions: It holds $4\leq v
\leq 64$, $24\leq v \leq 164$, $24\leq v \leq 264$ and $24\leq v
\leq 384$ when $d=15$, $d=22$, $d=29$ and $d=36$ respectively.
\end{prf}
We remark that all the Lie groups in the theorem have rank at least
$3$. Note that if one focusses on simple groups $G$, then only
dimensions $15$ and $36$ occur. Moreover, we easily derive the
following recognition theorem for the quaternionic volume:
\begin{cor}\label{INTcor02}
A $20$-dimensional Positive Quaternion K\"ahler Manifold with $\hat
A(M)[M]=0$ and quaternionic volume $v>384$ is symmetric.
\end{cor}
\vproof

\vspace{3mm}

The $\hat A$-genus of $M$ fits into the elliptic genus $\Phi(M)$ as
a first coefficient. The latter vanishes on the rationally
$3$-connected Wolf spaces of dimension $4n$ with odd $n$. This
implies $\hat A(M)[M]=\hat A(M,T_\cc M)[M]=\sign(M)=0$ on these
spaces. So it is reasonable to assume this last chain of equations.
This assumption can be combined with further slight conditions to
yield a confirmation of the main conjecture in dimension $20$. Note
that the original techniques in \cite{herr}---before we spotted the
error in there---would have permitted to derive this vanishing of
indices and to actually confirm the conjecture (in dimension $20$)
under a very mild upper bound on the Euler characteristic.


\section{Isometry groups}\label{sec05}

In this section we finish the proof of theorem \ref{theoA}. By
theorem \ref{INTtheo04} and table \ref{INTtable01} it remains to
show that a Positive Quaternion K\"ahler Manifold of dimension $20$
with isometry group one of $\Sp(4)$ or $\SO(9)$ up to finite
coverings is homothetic to the real Grassmannian $\widetilde
\Gr_4(\rr^9)$. The arguments given in this section will be Lie
theoretic mainly.

\vspace{3mm}

\textbf{In this section we obey the following convention: Every
statement on equalities, inclusions or decompositions of groups is a
statement up to finite coverings.}

\vspace{3mm}

First of all note the relation
\begin{align*}
\SO(n)\otimes \Sp(1)\cong (\SO(n)\times\Sp(1))/{\langle( (-1)^{n+1},
(-\id)^{n+1}) \rangle }
\end{align*}
This can be seen by computing the Lie algebras (cf.~\cite{koll},
p.~25) and determining the fibre of the covering $\SO(n)\times
\Sp(1)\to \SO(n)\otimes \Sp(1)$ over the identity.

From the tables in appendix B in \cite{koll}, p.~63--68, we cite
that all irreducible representations $\varrho(\G_2)$ of $\G_2$ in
degrees smaller than or equal to
\begin{itemize}
\item $\deg \varrho \leq 30$ if $\varrho\in \operatorname{Irr}_\rr(\G_2)$
\item $\deg \varrho \leq 15$ if $\varrho\in \operatorname{Irr}_\cc(\G_2)$
\item $\deg \varrho \leq 12$ if $\varrho\in \operatorname{Irr}_\hh(\G_2)$
\end{itemize}
are real and have degree
\begin{align}\label{CLAeqn01}
\deg \varrho \in \{7,14,27\}
\end{align}
This is considered an exemplary citation of the most important case.

We shall use this information in order to shed more light on
inclusions of Lie groups.
\begin{lemma}\label{CLAlemma01}
Let $G=G_1\times G_2$ be a decomposition of Lie groups. Let further
$H\neq 1$ be a simple Lie subgroup of $G$. Then (up to finite
coverings) $H$ is also a subgroup of one of $G_1$ and $G_2$ (by the
canonical projection).
\end{lemma}
\begin{prf}
Compose the inclusion $H\hto{} G_1\times G_2$ with the canonical
projection $G\to G_i$ (with $i\in \{1,2\}$) to obtain a morphism
$f_i: H \to G_i$. The kernel of $f_i$ is a normal subgroup of $H$,
i.e.~we have $\ker f_i\in\{1,H\}$ (up to finite coverings). If $\ker
f_i=1$, the morphism $f_i$ is an injection and we are done.
Otherwise, if $\ker f_i=H$, the morphism $f_i$ is constant. So if
both $f_1$ and $f_2$ are constant, the original inclusion $H\hto{}
G$ is a constant map, too. This contradicts $H\neq 1$.
\end{prf}

\begin{lemma}\label{CLAlemma02}
There is no inclusion of Lie groups $\SO(7)\hto{}\Sp(5)$ (not even
up to finite coverings).
\end{lemma}
\begin{prf}
By table \ref{PQKtable07} the group $\SO(7)$ has to be contained in
one of
\begin{align*}
\U(5),~ \Sp(4)\times \Sp(1),~ \Sp(3)\times \Sp(2),~ \SO(5)\otimes
\Sp(1)
\end{align*}
as there is no quaternionic representation of a simple Lie group $H$
in degree $10$ (other than the standard representation of $\Sp(5)$)
by the tables in appendix $B$ in \cite{koll}, p.~63--68. (The tables
neglect the cases of Lie groups with dimensions smaller than $11$.
Clearly, so can we.)

Thus, by lemma \ref{CLAlemma01} and for dimension reasons, we see
that $\SO(7)$ has to be a subgroup of one of
\begin{align*}
\SU(5),~ \Sp(4)
\end{align*}

Suppose first that $\SO(7)$ is a subgroup of $\SU(5)$. By table
\ref{PQKtable06}, lemma \ref{CLAlemma01} and for dimension reasons
this is not possible. (Again, there are no further irreducible
complex representations of degree $5$ of interest.)

Assume $\SO(7)$ is a subgroup of $\Sp(4)$. We argue in the analogous
way to get a contradiction. Alternatively, one may quote table
\ref{PQKtable09} for subgroups of maximal dimension.
\end{prf}
The following lemmas can fairly easily be proved using similar
arguments. We leave this to the reader.
\begin{lemma}\label{CLAlemma05}
There is no inclusion of Lie groups $\SU(5)\hto{}\SO(9)$ and
equally, $\SU(5)$ is not a subgroup of $\Sp(4)$ either---not even up
to finite coverings.
\end{lemma}
\vproof
\begin{lemma}\label{CLAlemma06}
There is no inclusion of Lie groups $\Sp(2)\times
\SU(3)\hto{}\SO(9)$ and equally, the group $\Sp(2)\times \SU(3)$
also is not a subgroup of $\Sp(4)$---not even up to finite
coverings.
\end{lemma}
\vproof
\newpage
\begin{lemma}\label{CLAlemma03}
Suppose that $k\in\{3,4\}$. The only inclusion of Lie groups
$\Sp(k)\hto{}\Sp(5)$ is given by the canonical blockwise inclusion
up to conjugation.
\end{lemma}
\vproof
\begin{lemma}\label{CLAlemma07}
There is no inclusion of Lie groups $\Sp(2)\times
\Sp(2)\hto{}\Sp(5)$ unless one of the $\Sp(2)$-factors includes
blockwise (up to conjugation).
\end{lemma}
\begin{prf}
By table \ref{PQKtable07} and by dimension $\Sp(2)\times \Sp(2)$
includes into one of
\begin{align*}
\U(5),~ \Sp(4)\times \Sp(1),~ \Sp(3)\times \Sp(2)
\end{align*}
Thus by lemma \ref{CLAlemma01} we see that $\Sp(2)\times \Sp(2)$ is
a subgroup of one of
\begin{align*}
\SU(5),~\Sp(3)\times \Sp(2)
\end{align*}
It cannot be a subgroup of $\SU(5)$, as a subgroup of $\SU(5)$ of
maximal dimension is of dimension $16<20=\dim\Sp(2)\times\Sp(2)$.
Thus $\Sp(2)\times \Sp(2)\In \Sp(3)\times \Sp(2)$ and we need to
determine all the possible inclusions $i$.

The only inclusion of $\Sp(2)$ into $\Sp(3)$ possible is given by
the canonical blockwise inclusion (up to conjugation). Suppose now
that both $\Sp(2)$-factors do not include blockwise. We shall lead
this to a contradiction.

We obtain that the inclusion composed with the canonical projection
\begin{align*}
i_1: \Sp(2)\times \{1\}\hto{} \Sp(3)\times \Sp(2)\to \Sp(2)
\end{align*}
again is an inclusion. This is due to the fact that the kernel of
this map has to be trivial, as $\Sp(2)$ is simple. Thus by our
assumption of non-blockwise inclusion the kernel has to be the
trivial group and $i_1$ is an inclusion. The same holds for
\begin{align*}
i_2: \{1\}\times \Sp(2)\hto{} \Sp(3)\times \Sp(2)\to\Sp(2)
\end{align*}
Without restriction, we may suppose that $i_1=i_2=\id$. Thus we
obtain that $i$ is an inclusion if and only if
\begin{align*}
\Sp(2)\times \Sp(2) \hto{} \Sp(3)\times \Sp(2)\to \Sp(3)
\end{align*}
is an inclusion. By consideration of rank this is impossible. This
yields a contradiction and at least one $\Sp(2)$-factor is
canonically included.
\end{prf}
\begin{lemma}\label{CLAlemma08}
The only inclusion of Lie groups $\SU(4)\hto{}\Sp(5)$ respectively
$\SU(4)\times \Sp(1)\hto{}\Sp(5)$ is given by the canonical
blockwise inclusion up to conjugation.
\end{lemma}
\begin{prf}
By table \ref{PQKtable07}, by the tables in appendix B in
\cite{koll}, p.~63--68, and by dimension we see that $\SU(4)$
respectively $\SU(4)\times \Sp(1)$ lies in one of
\begin{align*}
\U(5),~ \Sp(4)\times \Sp(1),~ \Sp(3)\times \Sp(2)
\end{align*}
The group $\SU(4)$ is not a subgroup of $\Sp(3)$ by table
\ref{PQKtable07} and by dimension. The only inclusion of $\SU(4)$
into $\U(5)$ is given by the canonical blockwise one due to the
usual arguments. The group $\SU(4)\times \Sp(1)$ is not a subgroup
of $\U(5)$; indeed, this is impossible by table \ref{PQKtable06},
the tables in appendix B in \cite{koll}, p.~63--68, and by
dimension.

\vspace{3mm}

\case{1} Thus in the case of the inclusion $\SU(4)\hto{}\Sp(5)$ we
observe that either the inclusion factors over $\SU(4)\hto{}
\U(4)\hto{} \Sp(5)$ and is given blockwise or the inclusion is given
via $\SU(4)\hto{}\Sp(4)\hto{}\Sp(5)$. Again, the inclusion
necessarily is given blockwise. For this we realise that $\SU(4)$
cannot be included into any maximal subgroup of $\Sp(4)$ other than
$\U(4)$.

\vspace{3mm}

\case{2} As for the inclusion $\SU(4)\times \Sp(1)\hto{}\Sp(5)$ we
note that $\SU(4)\times \Sp(1)$ maps into $\Sp(4)\times \Sp(1)$ by
dimension. By lemma \ref{CLAlemma01} and by dimension we see that
$\SU(4)$ again lies in $\Sp(4)$. As we have seen this inclusion
necessarily is blockwise. As $\Sp(4)$ maps into $\Sp(5)$ by
blockwise inclusion, the inclusion $\SU(4)\hto{}\Sp(5)$ is the
canonical one.

It then remains to see that $\Sp(1)$ includes into the
$\Sp(1)$-factor of $\Sp(4)\times \Sp(1)$ (and \emph{not} into the
$\Sp(4)$-factor). Assume this is not the case. This means that the
group $\Sp(1)$ necessarily does include into the $\Sp(4)$-factor.
Thus there is a homomorphism of groups
\begin{align*}
i: \SU(4)\times \Sp(1) \to \Sp(4)
\end{align*}
with the property that $i|_{\SU(4)}$ as well as $i|_{\Sp(1)}$ are
injective. (Clearly, the morphism $i$ itself cannot be injective.)
However, as we shall show, this contradicts the fact that $i$ is a
homomorphism: For $(x_1,x_2), (y_1,y_2)\in\SU(4)\times \Sp(1)$ we
compute
\begin{align*}
i((x_1,x_2)\cdot(y_1,y_2))&=i(x_1y_1,x_2y_2)=i(x_1,1)i(1,y_1)i(1,x_2)i(1,y_2)\\
i((x_1,x_2) \cdot
(y_1,y_2))&=i(x_1,x_2)i(y_1,y_2)=i(x_1,1)i(1,x_2)i(y_1,1)i(1,y_2)
\end{align*}
Thus we necessarily have that $i(1,y_1)i(1,x_2)=i(1,x_2)i(y_1,1)$.
As $i|_{\Sp(1)}$ is injective, we realise that---whatever the
inclusion $i|_{\Sp(1)}$ will be---the group $i(\Sp(1))$ always
contains elements that do not commute with every element of
$\SU(4)\In\Sp(4)$. (Clearly, $\Sp(1)$ cannot be included into the
centre $T^4$ of $\SU(4)$.)

Consequently, we obtain that the inclusion of $\Sp(1)$  into
$\Sp(4)\times \Sp(1)$ maps $\Sp(1)$ to the $\Sp(1)$-factor and the
projection
\begin{align*}
\Sp(1)\to \Sp(4)\times \Sp(1)\to\Sp(4)
\end{align*}
is the trivial
map. This proves the assertion.
\end{prf}

\vspace{3mm}

Let us now prove the classification result.
\begin{theo}\label{CLAtheo01}
A $20$-dimensional Positive Quaternion K\"ahler Manifold $M$ which
has an isometry group $\Isom(M)$ that satisfies
\begin{align*}
\Isom_0(M)\in \{\SO(9),\Sp(4)\}
\end{align*}
up to finite coverings is homothetic to the real Grassmannian
\begin{align*}
M\cong \widetilde \Gr_4(\rr^9)
\end{align*}
\end{theo}
\begin{prf}
We proceed in three steps. First we shall establish a list of
stabiliser groups in a $T^4$-fixed-point---where $T^4$ is the
maximal torus of $\Isom_0(M)$---that might occur unless $M$ is a
Wolf space. As a second step we reduce the list by inclusions into
the isometry group and the holonomy group. Finally, in the third
step we show by more distinguished arguments that also the remaining
stabilisers from the list cannot occur, whence $M$ has to be
symmetric.

\vspace{3mm}

\step{1} Both groups $\SO(9)$ as well as $\Sp(4)$ have rank $4$,
i.e.~they contain a $4$-torus $T^4=\s^1\times\dots\times \s^1$. The
Positive Quaternion K\"ahler Manifold $M$ has positive Euler
characteristic by theorem \ref{PQKtheo05}. Thus by the Lefschetz
fixed-point theorem we derive that there exists a $T^4$-fixed-point
$x\in M$. Let $H_x$ denote the (identity-component of the) isotropy
group of the $G$-action in $x$ for $G\in\{\SO(9),\Sp(4)\}$. Since
$G$ is of dimension $36$ and since $\dim M=20$, we obtain that $\dim
H_x\geq 17$ unless the action of $G$ on $M$ is transitive. If this
is the case, a result by Alekseevskii (cf.~\cite{bess}.14.56,
p.~409) yields the symmetry of $M$. If $M$ is symmetric, it is a
Wolf space and the dimension of the isometry group then yields that
$M\cong \widetilde \Gr_4(\rr^9)$. We may even assume $\dim H_x\geq
18$ by the classification of cohomogeneity one Positive Quaternion
K\"ahler Manifolds in theorem \cite{danc2}.7.4, p.~24.

Since $\rk G=4$ and since $x$ is a $T^4$-fixed-point, we also obtain
$\rk H_x=4$. Moreover, $H_x$ is a closed subgroup. Thus we may give
a list of all products of semi-simple Lie groups and tori of
dimension $36\geq \dim H_x\geq 18$ and with $\rk H_x=4$ up to finite
coverings:
\begin{align*}
&\SU(5),~ \SO(9),~ \SO(8),~ \Sp(4),~ \Sp(3)\times \Sp(1),~
\Sp(3)\times \s^1,\\& \SO(7)\times \Sp(1),~ \SO(7)\times \s^1, ~
\SO(6)\times \Sp(1),~ \Sp(2)\times \Sp(2),\\& \Sp(2)\times \SU(3),~
\Sp(2)\times \G_2,~ \G_2\times \G_2,~ \G_2\times \SU(3),\\& \G_2\times \Sp(1)\times
\Sp(1),~ \G_2\times \Sp(1)\times \s^1
\end{align*}

\vspace{3mm}

\step{2} We now apply two criteria by which we may reduce the list:
\begin{itemize}
\item
On the one hand we have that $H_x$ is a Lie subgroup of $G$.
\item
On the other hand by the isotropy representation $H_x$ is a Lie
subgroup of $\Sp(5)\Sp(1)$---cf.~theorem \cite{kob1}.VI.4.6, p.~248.
\end{itemize}
We use lemma \ref{CLAlemma01} to see that every group $H_x$ in the
list contains a factor that has to include into $\Sp(5)$ up to
finite coverings.

An iterative application of table \ref{PQKtable08} yields that every
maximal rank subgroup of the classical group $G\in
\{\SO(9),\Sp(4)\}$ again is a product of classical groups. Thus
$H_x$ may not be one of the groups
\begin{align*}
&\Sp(2)\times \G_2,~
 \G_2\times \G_2,~ \G_2\times \SU(3),\\& \G_2\times \Sp(1)\times
\Sp(1),~ \G_2\times \Sp(1)\times \s^1
\end{align*}
(not even up to finite coverings).

Now apply lemma \ref{CLAlemma02} in the respective cases to reduce
the list of potential stabilisers to
\begin{align*}
&\SU(5),~ \Sp(4),~ \Sp(3)\times \Sp(1),~ \Sp(3)\times \s^1,\\&
 \SO(6)\times \Sp(1),~
\Sp(2)\times \Sp(2),~ \Sp(2)\times \SU(3)
\end{align*}
Indeed, this lemma rules out all the groups $H_x$ that contain a
factor of the form $\SO(7)$; and as we see that $\SO(7)$ is not a
subgroup of $\Sp(5)$, also $\SO(8)$ and $\SO(9)$ cannot be subgroups
of $\Sp(5)$.

Now apply lemmas \ref{CLAlemma05} and \ref{CLAlemma06} by which
potential inclusions into the isometry group $G$ are made clearer.
That is, they rule out the groups $\SU(5)$ and $\Sp(2)\times
\SU(3)$. Thus the list of possible isotropy groups reduces further
to
\begin{align*}
\Sp(4),~ \Sp(3)\times \Sp(1),~ \Sp(3)\times \s^1,~
 \SO(6)\times \Sp(1),~
\Sp(2)\times \Sp(2)
\end{align*}

\vspace{3mm}

\step{3} Let us consider the inclusion of the groups $H_x$ into the
holonomy group $\Sp(5)\Sp(1)$. By \ref{CLAlemma01} the largest
direct factor of the candidates in our list has to be a subgroup of
$\Sp(5)$ (up to finite coverings), as it cannot be included into
$\Sp(1)$. By lemma \ref{CLAlemma03} we see that the inclusion of
$\Sp(4)$ into $\Sp(5)$ and the one of the $\Sp(3)$-factor of
$\Sp(3)\times \Sp(1)$ respectively of $\Sp(3)\times \s^{1}$ has to
be blockwise. Due to lemma \ref{CLAlemma07} we observe that there is
also an $\Sp(2)$-factor of $\Sp(2)\times \Sp(2)$ that includes
blockwise into $\Sp(5)$. Thus we obtain that every group from our
list which contains a factor of the form $\Sp(k)$ for $k\geq 2$ has
a circle subgroup $\s^1\In \Sp(k)$ that includes into $\Sp(5)$ by
\begin{align*}
\diag(\s^1,1,1,1,1) \In \diag(\Sp(k),1,\dots,1) \In \Sp(5)
\end{align*}
Thus this circle group fixes a codimension $4$ Positive Quaternion
K\"ahler component. Due to theorem \cite{fang2}.1.2, p.~2, we obtain
that $M\cong \hh\pp^5$ or $M\cong \Gr_2(\cc^{7})$; a contradiction
by our assumption on the isometry group.

\vspace{3mm}

This leaves us with $H_x=\SO(6)\times \SO(3)$, which is
$\SU(4)\times \Sp(1)$ up to $(\zz_2\oplus\zz_2)$-covering.
Equivalently, we consider an orbit of the form
\begin{align*}
X:&= \frac{\SO(9)}{\SO(6)\times \SO(3)}\\
&=\frac{\Spin(9)}{(\Spin(6)\times \Spin(3))/\langle
(-\id,-\id)\rangle} \\&= \frac{\Spin(9)}{(\SU(4)\times
\Sp(1))/\langle (-\id,-1)\rangle}
\end{align*}
Now consider the isotropy representation of the stabiliser group.
There are basically two possibilities: Either the whole stabiliser
includes into the $\Sp(5)$-factor or the $\Sp(1)$-factor includes
into the $\Sp(1)$-factor of $\Sp(5)\Sp(1)$ (and not into the
$\Sp(5)$-factor).

In the first case we apply lemma \ref{CLAlemma08} to see that the
inclusion of $\SU(4)\times \Sp(1)$ into $\Sp(5)$ is blockwise. So is
the inclusion of the $\Sp(1)$-factor in particular. Thus again we
obtain a sphere which is represented by
\linebreak[4]$\diag(\s^1,1,1,1,1)$ and which fixes a codimension
four quaternionic component. We proceed as above.

Let us now deal with the second case. Again we cite lemma
\ref{CLAlemma08} to see that the $\SU(4)$-factor includes into
$\Sp(5)$ in a blockwise way. Observe now that the tangent bundle
$TM$ of $M$ splits as
\begin{align*}
TM=TX\oplus NX
\end{align*}
over $X$, where $NX$ denotes the normal bundle.
 Since
\begin{align*}
\dim X=\dim \SO(9)-\dim \SO(6)\times\SO(3)=36-18=18
\end{align*}
we obtain that the normal bundle is two-dimensional. Thus the slice
representation of the isotropy group $(\SU(4)\times
\Sp(1))/(-\id,-1)$ at a fixed-point, i.e.~the representation on
$NX$, is necessarily trivial. That is, the action of the isotropy
group $(\SU(4)\times \Sp(1))/(-\id,-1)$ at a fixed-point has to
leave the normal bundle pointwise fixed.

The $\Sp(1)$-factor of $\SU(4)\times \Sp(1)$ maps isomorphically (up
to finite coverings) into the $\Sp(1)$-factor of the holonomy group;
the $\SU(4)$-factor maps into $\Sp(5)$. Thus the action of this
$\SU(4)\times \Sp(1)$ on the tangent space $T_xM\cong \hh^5$ at
$x\in M$ is given by $(A,h)(v)=Avh^{-1}$. This action, however, has
no $18$-dimensional (respectively $2$-dimensional) invariant
subspace as the $\Sp(1)$-factor acts transitively on each
$\hh$-component. Hence the normal bundle does not remain fixed under
the action of the stabiliser. Thus $\SU(4)\times \Sp(1)$ cannot
occur as an isotropy group.

\vspace{3mm}

Hence we have excluded all the cases that arose from the assumption
$\dim H_x\geq 18$. This was equivalent to the action of the isometry
group neither being transitive nor of cohomogeneity one. In the
latter two cases---as already observed---the manifold $M$ has to be
symmetric. More precisely, since there is no $20$-dimensional Wolf
space with $\Isom_0(M)=\Sp(4)$ (up to finite coverings), we obtain
$\Isom_0 M=\SO(9)$ and $M\cong \widetilde \Gr_4(\rr^9)$.
\end{prf}
Theorem \ref{theoA} now follows from theorems \ref{INTtheo04} and
\ref{CLAtheo01}.

\vspace{3mm}

Clearly, as for $\dim \Isom_0(M)\not\in \{15,22,29\}$ one hopes an
approach by similar techniques as in the proof of \ref{CLAtheo01} to
be likewise successful. Yet, we remark that for example in dimension
$29$ one will have to cope with five different isometry groups due
to table \ref{INTtable01}. All these groups are of rank $5$. So one
lists all the possible stabilisers at a $T^5$-fixed-point on $M$
that do not necessarily make the action of $\Isom(M)$ transitive or
of cohomogeneity one. That is, one computes all the products $H$ of
semi-simple Lie groups and tori that satisfy $\rk H=5$ and $11\leq
\dim H \leq 29$. This results in a list of $45$ possible groups $H$
(up to finite coverings). Following our previous line of argument we
then try to rule out stabilisers by showing that they either may not
include into a respective isometry group or that they may not be a
subgroup of $\Sp(5)\Sp(1)$. If both is not the case, as a next step
we try to show that the way $H$ includes into the holonomy group
already implies the existence of an $\s^1$-fixed-point component of
codimension $4$. This would imply the symmetry of the ambient
manifold $M^{20}$. We observe that by far the biggest part of this
procedure is covered by the arguments we applied before and we
encourage the reader to provide the concrete reasoning. Nonetheless,
we encounter new difficulties: For example, the group $\SU(4)\times
\s^1\times \s^1$ includes into the isometry group $\SU(4)\times
\G_2\cong \SO(6)\times \G_2$. If its inclusion into $\Sp(5)\Sp(1)$
is induced by the blockwise inclusion of $\SU(4)$ into $\Sp(5)$, the
canonical inclusion of $\s^1$ into $\Sp(1)$ and the diagonal
inclusion of $\s^1$ into $\Sp(5)$, we realise that there is no
codimension four $\s^1$-fixed-point component. Then methods more
particular in nature will have to be provided---as we did in step 3
of the proof of theorem \ref{CLAtheo01}. We leave this to the
reader.


\section{Proof of theorem C}\label{sec06}

Since one may relate the index $i^{0,n+2}$ (cf.~\ref{PQKtheo10})
directly to the dimension of $\Isom(M^{4n})$, it seems to be pretty
natural to try to provide a recognition theorem on this information.
This will result in theorem \ref{theoC}, the first one to identify
the real Grassmanian.

For the proof of theorem \ref{theoC} we imitate and generalise the
techniques and results from section \ref{sec05}. \textbf{Again, we
shall neglect finite coverings.}
\begin{lemma}\label{REClemma01}
For $n\geq 6$ there is no inclusion of Lie groups
$\SO(n+1)\hto{}\Sp(n)$, not even up to finite coverings.
\end{lemma}
\begin{prf}
Due to table \ref{PQKtable07} the group $\SO(n+1)$ either has to be
contained in $\U(n)$, $\Sp(k)\times \Sp(n-k)$ with $1\leq k\leq
n-1$, some $\SO(p)\otimes \Sp(q)$ with $pq=n$, $p\geq 3$, $q\geq 1$
or in $\varrho(H)$ for a simple Lie group $H$ and an irreducible
quaternionic representation $\varrho\in \operatorname{Irr}_\hh(H)$
of dimension $\deg \varrho=2n$. The cases with direct product or
tensor product yield an inclusion of $\SO(n+1)$ in either some
$\SO(k)$ with $k\geq n$---which is impossible by dimension---or into
some smaller symplectic group by lemma \ref{CLAlemma01}.

Assume there is an inclusion into $\U(n)=(\SU(n)\times
\U(1))/\zz_n$. Then again lemma \ref{CLAlemma01} yields an inclusion
into $\SU(n)$. By table \ref{PQKtable06} the maximal subgroups of
$\SU(n)$ are given by $\SO(n)$, $\Sp(m)$ with $2m=n$,
$\mathbf{S}(\U(k)\times \U(n-k))$ for ($1\leq k \leq n-1$),
$\SU(p)\otimes \SU(q)$ with $pq=n$, $p\geq 3$, $q\geq 2$ and by
$\varrho(H)$ for a simple Lie group $H$ and an irreducible
quaternionic representation $\varrho\in \operatorname{Irr}_\cc(H)$
of dimension $\deg \varrho=n$. An inclusion in the first case is
impossible due to dimension. Cases two to four lead to inclusions
into smaller symplectic or special unitary groups by lemma
\ref{CLAlemma01}.

\vspace{3mm}

Hence we need to have a closer look at irreducible quaternionic and
complex representations of simple Lie groups $H$. The tables in
\cite{koll}, appendix B, p.~63--68, give all the representations of
simple Lie groups satisfying a certain dimension bound, which is
given by
\begin{align*}
2 \dim H&\geq \deg \varrho -2\\
\dim H&\geq \deg \varrho-1\\
\dim H&\geq \frac{3}{2}\deg \varrho-4
\end{align*}
for real, complex and quaternionic representations respectively.

First of all for $n\geq 3$ the tables together with our previous
reasoning yield that $k=n$ is the maximal number for which $\SU(k)$
is a maximal subgroup of $\Sp(n)$. Equally, for $n\geq 7$ we obtain
that $k=n$ is the maximal number for which $\SO(k)$ is a maximal
subgroup of $\Sp(n)$ or of $\SU(n)$. This means in particular that
$\SO(n+1)$ cannot be included into $\Sp(n)$ by a chain
\begin{align}\label{RECeqn02}
\SO(n+1) \In G_1 \In  \dots \In G_l \In \Sp(n)
\end{align}
of (irreducible representations of) classical groups $G_1,\dots,
G_l$ for $n\geq 6$.

It remains to prove that there is no such chain involving
(representations of) exceptional Lie groups $G_i$. For this it
suffices to realise that there are no exceptional Lie groups $H$
satisfying
\begin{align}\label{RECeqn01}
\dim \SO(n+1) \leq \dim H\leq \dim \Sp(n)
\end{align}
with $H$ admitting a quaternionic or complex representation of
degree smaller than or equal to $2n$ or $n$ respectively. (We
clearly may neglect the real representations $\varrho$ of degree $k$
with $k\leq n$, as there evidently cannot be inclusions $\SO(n+1)\In
\varrho(H)\In \SO(k)$ by dimension.)

In table \ref{RECtable01} for each exceptional Lie group $H$ we give
the values of $n$ for which the inequalities \eqref{RECeqn01} are
satisfied. Additionally, we note the corresponding maximal degree
$\deg \varrho=2n$ ($\deg \varrho=n$) of an irreducible quaternionic
(complex) representation $\varrho$  by which $H$ might become the
subgroup $\varrho(H)\In \Sp(k)$ ($\varrho(H)\In\SU(k)$) with $k\leq
n$ for the given values of $n$. That is, for example in the case of
$\G_2$ we see that if there is a quaternionic representation (a
complex representation) of degree smaller than or equal to $22$ (to
$11$), then there is an inclusion of $\G_2$ into $\Sp(k)$ (into
$\SU(k)$) for $k\leq 11$ and now also conversely: If there is no
such representation, then $\G_2$ cannot be a subgroup  satisfying
$\SO(n+1)\In \G_2 \In \Sp(n)$ for any $n\in \nn$.

Now the tables in \cite{koll} yield that there are no quaternionic
respectively complex representations of $H$ in the degrees depicted
in table \ref{RECtable01}. This amounts to the fact that for the
relevant values of $n$ from table \ref{RECtable01} there is no
inclusion $\SO(n+1)\In H \In \Sp(n)$. Thus by \eqref{RECeqn01} there
are no inclusions of exceptional Lie groups $H$ with $\SO(n+1)\In
H\In \Sp(n)$ for any $n\in \nn$.

Thus we have proved that there cannot be a chain of the form
\eqref{RECeqn02} with an exceptional Lie group $G_i$. Combining this
with our previous arguments proves the assertion.
\begin{table}[h]
\centering \caption{Degrees of relevant representations}
\label{RECtable01} \begin{tabular}{|@{\hspace{6mm}}c@{\hspace{6mm}}|
@{\hspace{6mm}}c@{\hspace{6mm}}|@{\hspace{6mm}}c@{\hspace{6mm}}| }
\hline Lie group & $n\in$ & $\deg \varrho \leq$  \\
\hline\hline $\G_2$ & $\{3,4\}$ & $8, 4$
\\ \hline
$\F_4$& $\{5,6,7,8,9\}$& $18, 9$ \\
\hline $\E_6$ & $\{7,8,9,10,11\}$ & $22, 11$ \\
\hline $\E_7$ & $\{8,9,10,11,12,13,14,15\}$ & $30, 15$\\\hline
$\E_8$ & $\{11,12,\dots ,20,21\}$ & $42, 21$
\\\hline
\end{tabular}
\end{table}
\end{prf}
Note that the bound $n\geq 7$ in the lemma is necessary since the
universal two-sheeted covering of $\SO(6)$ is $\SU(4)$. In higher
dimensions no such exceptional identities occur as can be seen from
the corresponding Dynkin diagrams.
\begin{lemma}\label{REClemma02}
For $n\geq 3$ the only inclusion of Lie groups $\Sp(\lfloor
\frac{n}{2}\rfloor+1)\hto{}\Sp(n)$ is given by the canonical
blockwise one up to conjugation.
\end{lemma}
\begin{prf}
We proceed as in lemma \ref{REClemma01}. Indeed, by the same
arguments as above we see that every chain of classical groups
\begin{align}\label{RECeqn03}
\Sp\bigg(\bigg\lfloor \frac{n}{2}\bigg\rfloor+1\bigg) \In G_1 \In
\dots \In G_l \In\Sp(n)
\end{align}
involves symplectic or special unitary groups $G_i$ of rank smaller
than or equal to $n$ only. For this we use that there is no
inclusion $\Sp(\lfloor \frac{n}{2}\rfloor+1)\In \SO(n)$ by
dimension; indeed
\begin{align*}
\dim \SO(n)=\frac{n(n-1)}{2}&< \left\{\begin{array}{cl}
\big(\frac{n}{2}+1\big)(n+3) & \textrm{for } n \textrm{ even}\\
\frac{(n+1)(n+2)}{2} & \textrm{for } n \textrm{ odd}
\end{array}
\right\}\\&=\dim \Sp\bigg(\bigg\lfloor \frac{n}{2}\bigg\rfloor+1\bigg)
\end{align*}
More precisely, in such a chain of classical groups the inclusion of
$\Sp(\lfloor \frac{n}{2}\rfloor+1)$ is necessarily blockwise since
$n\geq 3$. This is due to the fact that actually only symplectic
groups $G_i$ which are included in a blockwise way may appear by
table \ref{PQKtable06}; i.e.~the subgroups of $\SU(n)$ are to small
to permit the inclusion of $\Sp(\lfloor \frac{n}{2}\rfloor+1)$.

We now have to realise that there is no chain as in \eqref{RECeqn03}
with an exceptional Lie group $G_i$. As in the proof of lemma
\ref{REClemma01} we depict the values of $n$ for which an inclusion
$\Sp(\lfloor \frac{n}{2}\rfloor+1) \In H\In \Sp(n)$ of an
exceptional Lie group $H$ might be possible---when merely
considering dimensions---in table \ref{RECtable02}. The table again
also yields degree bounds for the degrees of quaternionic and
complex representations.
\begin{table}[h]
\centering \caption{Degrees of relevant representations}
\label{RECtable02}
\begin{tabular}{|@{\hspace{6mm}}c@{\hspace{6mm}}| @{\hspace{6mm}}c@{\hspace{6mm}}|@{\hspace{6mm}}c@{\hspace{6mm}}| }
\hline Lie group & $n\in$ & $\deg \varrho \leq$  \\
\hline\hline $\G_2$ & $\{3\}$ & $6, 3$
\\ \hline
$\F_4$& $\{5,6,7\}$& $14, 7$ \\
\hline $\E_6$ & $\{7,8,9\}$ & $18, 9$ \\
\hline $\E_7$ & $\{8,9,10,11,12,13\}$ & $26, 13$\\\hline $\E_8$ &
$\{11,12, 13, 14, 15, 16, 17, 18, 19\}$ & $38, 19$
\\\hline
\end{tabular}
\end{table}
 Then the tables in appendix \cite{koll}.B,
p.~63--68, yield that under these respective restrictions no
representations of exceptional Lie groups can be found. This implies
that each $G_i$ in the chain is classical. Thus the inclusion of
$\Sp(\lfloor \frac{n}{2}\rfloor+1)$ into $\Sp(n)$ is necessarily
given blockwise.
\end{prf}
Note that $\Sp(1)\cong \SU(2)$, whence the inclusion
$\Sp(1)\hto{}\Sp(2)$ is not necessarily blockwise.
\begin{lemma}\label{REClemma04}
Let $n\geq 7$ be odd. Every inclusion $\Sp(\frac{n+1}{2}-1)\times
\Sp(2)\hto{}\Sp(n)$ restricts to the canonical blockwise one (up to
conjugation and finite coverings) on the first factor.
\end{lemma}
\begin{prf}
By table \ref{PQKtable07} and by dimension we see that
$\Sp(\frac{n+1}{2}-1)\times \Sp(2)$ lies in one of
\begin{align*}
\U(n),~ \Sp(k)\times \Sp(n-k) \textrm{ for } 1\leq k\leq n-1,~
\varrho(H)
\end{align*}
for a simple Lie group $H$ and an irreducible quaternionic
representation $\varrho$ of degree $\deg \varrho=2n$.

If $\Sp(\frac{n+1}{2}-1)\times \Sp(2)$ should happen to appear as a
subgroup of $\SU(n)$, then table \ref{PQKtable06} would show that
\begin{align*}
\Sp\bigg(\frac{n+1}{2}-1\bigg)\hto{} \SU(n-1) \hto{} \SU(n)
\end{align*}
necessarily is included in the standard ``diagonal'' way induced by
the standard inclusion $\hh\hto{}\cc^{2\times 2}$. For this we
observe the following facts that result when additionally taking
into account the tables in appendix \cite{koll}.B, p.~63--68: The
largest special orthogonal subgroup (up to finite coverings) of
$\SU(n)$ is $\SO(n)$ for $n\geq 6$. The group $\SO(n)$ does not
permit $\Sp(\frac{n+1}{2}-1)$ as a subgroup for $n\geq 4$. Moreover,
there are no irreducible complex representations by which a simple
Lie group  $H$ might include into some $\SU(k)$ (for $k\leq n$)
satisfying $\Sp(\frac{n+1}{2}-1)\In H$.

Now we see that whenever $\Sp(\frac{n+1}{2}-1)$ is included
diagonally into $\SU(n)$ as depicted, there is no inclusion of
$\Sp(\frac{n+1}{2}-1)\times \Sp(2)$ possible. That is, for the
inclusion of this direct product to be a homomorphism we need the
group $\Sp(2)$ to map into the centraliser
$C_{\SU(n)}(\Sp(\frac{n+1}{2}-1))$ of $\Sp(\frac{n+1}{2}-1)$ in
$\SU(n)$. Yet, we obtain
\begin{align*}
C_{\SU(n)}\bigg(\Sp\bigg(\frac{n+1}{2}-1\bigg)\bigg)\cong \s^1\times
\s^1
\end{align*}
and thus no inclusion of $\Sp(2)$ is possible. Hence
$\Sp(\frac{n+1}{2}-1)\times \Sp(2)$ cannot be a subgroup of $\U(n)$.

\vspace{3mm}

The tables in \cite{koll} again yield that whenever a simple Lie
group $H$ is included into $\Sp(k)$ (for $k\leq n$) via an
irreducible quaternionic representation $\varrho$, the inclusion is
one of
\begin{align*}
\SU(6)\hto{} \Sp(10),~ \SO(11)\hto{}\Sp(16),~ \SO(12)\hto{}\Sp(16),~
 \E_7\hto{} \Sp(28)
\end{align*}
or an inclusion of a symplectic group of rank smaller than $k$
unless the degree of the representation is far too large to be of
interest for our purposes. Indeed, already the depicted inclusions
are not relevant, since $\Sp(\frac{n+1}{2}-1)$ cannot be included
into $\SU(6),~ \SO(11),~ \SO(12),~ \E_7$ respectively when $n\geq
10,16,16,28$.

\vspace{3mm}

Thus we see that every inclusion of $\Sp(\frac{n+1}{2}-1)\times
\Sp(2)$ has to factor through one of $\Sp(n-k)\times \Sp(k)$ for
$1\leq k\leq n-1$. Hence for $1\leq k<\frac{n+1}{2}-1$ we obtain
that the only inclusion of $\Sp(\frac{n+1}{2}-1)$ into $\Sp(n)$
factoring through $\Sp(n-k)\times \Sp(k)$ is given by the standard
blockwise inclusion
\begin{align*} \Sp\bigg(\frac{n+1}{2}-1\bigg)\hto{}
\Sp(n-k)\hto{}\Sp(k)
\end{align*}
Thus we may assume without restriction that $k=\frac{n+1}{2}-1$ (and
$n-k=\frac{n+1}{2}$) and that the inclusion of
$\Sp(\frac{n+1}{2}-1)\hto{}\Sp(n)$ is not the standard blockwise
one. Thus we see that the inclusion necessarily factors over
\begin{align*}
\Sp\bigg(\frac{n+1}{2}-1\bigg)\hto{}\Sp\bigg(\frac{n+1}{2}-1\bigg)\times\Sp\bigg(\frac{n+1}{2}\bigg)
\hto{}\Sp(n)
\end{align*}
where the first inclusion splits as a product of the standard
blockwise inclusions
\begin{align*}
\Sp\bigg(\frac{n+1}{2}-1\bigg)&\hto{\id}\Sp\bigg(\frac{n+1}{2}-1\bigg)\\
\Sp\bigg(\frac{n+1}{2}-1\bigg)&\hto{} \Sp\bigg(\frac{n+1}{2}\bigg)
\end{align*}
So regard $\Sp(\frac{n+1}{2}-1)$ as the subgroup of $\Sp(n)$ given
by this inclusion. Again we make use of the fact that the
$\Sp(2)$-factor has to include into the centraliser
\begin{align*}
C_{\Sp(n)}\bigg(\Sp\bigg(\frac{n+1}{2}-1\bigg)\bigg)\cong \s^1\times
\s^1\times \Sp(1)
\end{align*}
Such an inclusion clearly is impossible and we obtain a
contradiction. Thus the inclusion of $\Sp(\frac{n+1}{2}-1)\times
\Sp(2)$ is the standard blockwise one when restricted to the
$\Sp(\frac{n+1}{2}-1)$-factor.
\end{prf}
\begin{lemma}\label{REClemma03}
In table \ref{RECtable04} the semi-simple Lie groups of maximal
dimension with respect to a fixed rank (from rank $1$ to rank $12$)
are given up to isomorphisms and finite coverings. From rank $13$ on
the groups that are maximal in this sense are given by the two
infinite series $\Sp(n)$ and $\SO(2n+1)$ only.
\begin{table}[h]
\centering \caption{Largest Lie groups with respect to fixed rank}
\label{RECtable04}
\begin{tabular}{|@{\hspace{6mm}}r@{\hspace{6mm}}| @{\hspace{6mm}}c@{\hspace{6mm}}|@{\hspace{6mm}}r@{\hspace{6mm}}| }
\hline $\rk G$ & extremal Lie groups $G$ & $\dim G$  \\
\hline\hline $1$ & $\Sp(1)$ & $3$
\\ \hline
$2$& $\G_2$& $14$ \\
\hline $3$ & $\Sp(3), \SO(7)$ & $21$ \\
\hline $4$ & $\F_4$ & $52$\\\hline $5$ & $\Sp(5), \F_4\times
\Sp(1),\SO(11)$ & $55$
\\\hline
 $6$ & $\E_6,\Sp(6)$ & $78$\\
\hline $7$ & $\E_7$ & $133$ \\\hline $8$ & $\E_8$ & $248$\\
\hline $9$ & $\E_8\times \Sp(1)$ & $251$\\ \hline $10$ & $\E_8\times
\G_2$ & $262$\\ \hline $11$ & $\E_8\times \Sp(3), \E_8\times \SO(7)$
& $269$\\\hline $12$ & $\Sp(12), \SO(25), \E_8\times \F_4$ &
$300$\\\hline
\end{tabular}
\end{table}
\end{lemma}
\begin{prf}
Table \ref{RECtable04} results from a case by case check using table
\ref{PQKtable04}. From dimension $13$ on semi-simple Lie groups
involving factors that are exceptional Lie groups are smaller than
the largest classical groups. Among products of classical groups the
ratio between dimension and rank is maximal for the types $\B_n$ and
$\C_n$,
\end{prf}

Now we provided all the tools that will permit to prove theorem \ref{theoC}.
\begin{proof}[\textsc{Proof of theorem \ref{theoC}}]
\step{1} By theorem \ref{PQKtheo07} we may suppose that $\rk
\Isom(M^{4n})\leq \lceil \frac{n}{2}\rceil+2$, since otherwise $M\in
\{\hh\pp^n,\Gr_2(\cc^{n+2})\}$. The dimension bounds in table
\ref{tableC} for $4\leq n\leq 20$ result from table
\ref{RECtable04}: That is, for each such $n$ the bound is the
dimension of the largest group that satisfies this rank condition.
Thus every group with larger dimension has rank large enough to
identify $M^{4n}$ as one of $\hh\pp^n$ and $\Gr_2(\cc^{n+2})$.

In degree $n=3$ we use that there is no semi-simple Lie group of
rank smaller than or equal to $4$ in dimensions $29$ to $36=\dim
\Sp(4)$. By theorem \ref{PQKtheo07} we have $\dim \Isom(M^{12})\leq
36$.

\vspace{3mm}

\step{2} Now we determine all (the one-components of) the isometry
groups $G=\Isom_0(M^{4n})$ (up to finite coverings) with $\rk G\leq
\lceil \frac{n}{2}\rceil+2$ satisfying the dimension bound for
$n\geq 22$ and  $n\not\in\{27,28\}$ as
\begin{align*}
G\in \bigg\{&\SO(n+4),~ \SO(n+5),~ \Sp\bigg(\frac{n}{2}+2\bigg)
\bigg\} \intertext{for $n$  even and as} G\in\bigg\{&\SO(n+4),~
\SO(n+4)\times \SO(2),~ \SO(n+4)\times \SO(3),\\&
\SO(n+5),~\SO(n+6),~
\Sp\bigg(\frac{n+1}{2}+2\bigg),~\Sp\bigg(\frac{n+1}{2}+1\bigg),\\& \Sp\bigg(\frac{n+1}{2}+1\bigg)\times
\SO(2),~ \Sp\bigg(\frac{n+1}{2}+1\bigg)\times \Sp(1)\bigg\}
\end{align*}
for $n$ odd. This can be achieved as follows: We see that whenever
we have a product of classical groups we may replace it by a simple
classical group of the same rank and of larger dimension. The
classical groups for which the ratio between dimension and rank is
maximal are given by the groups of type $\B$ and $\C$. Moreover, the
series $\dim \B_n=\dim \C_n$ is strictly increasing in $n$. We
compute
\begin{align*}
\dim \SO(n+3)\times \SO(3)&=\frac{n^2+5n+12}{2} \intertext{whilst}
\rk \SO(n+3)\times \SO(3)&=\bigg\lfloor
\frac{n+3}{2}\bigg\rfloor+1=\bigg\lceil \frac{n}{2}\bigg\rceil+2
\end{align*}
Consequently, by our assumption on the dimension of $G$ we need to
find all the groups that are larger in dimension but not larger in
rank than $\SO(n+3)\times \SO(3)$. This process results in the list
we gave.

We still need to see when there are groups $G$ that are larger in
dimension than $\SO(n+3)\times \SO(3)$ but not larger in rank and
that have exceptional Lie groups as direct factors (up to finite
coverings). Clearly, this can only happen in low dimensions. So we
use lemma \ref{REClemma03} and table \ref{RECtable04} to see that
unless $n\in\{27,28\}$ there do no appear exceptional Lie groups as
factors. As for degrees $n\in \{27,28\}$ the group $\E_8\times \E_8$
has dimension
\begin{align*}
\dim (\E_8\times \E_8)=496>\begin{cases} 438=\dim \SO(30)\times
\SO(3)\\ 468=\dim \SO(31)\times \SO(3)
\end{cases}
\end{align*}
Therefore in these degrees we want to assume that $\dim
\Isom(M^{4n})>496$. This will make it impossible to identify the
real Grassmannian, since \linebreak[4]$\dim \SO(31)=465$ and since $\dim
\SO(32)=496$. Nonetheless the following arguments hold as well.

\vspace{3mm}

\step{3} We now prove that whenever $G$ is taken out of the list we
gave, then actually $G=\SO(n+4)$ and $M\cong \widetilde
\Gr_4(\rr^{n+4})$. In order to establish this we shall have a closer
look at orbits around a fixed-point of the maximal torus of $G$ for
each respective possibility of $G$. (Such a point exists due to the
Lefschetz fixed-point theorem and the fact that
$\chi(M)>0$---cf.~\ref{PQKtheo05}.) This will lead to the
observation that a potential action of $G$ has to be transitive,
whence $M$ is homogeneous. Due to Alekseevski homogeneous Positive
Quaternion K\"ahler Manifolds are Wolf
spaces---cf.~\cite{bess}.14.56, p.~409.

Since $\dim M=4n$, we necessarily obtain that the orbit $G/H$ of $G$
has dimension at most $\dim G/H\leq 4n$. Assume first that $G$ is a
direct product from the list with $\SO(n+4)$ as a factor (up to
finite coverings). Thus all the maximal rank subgroups $H$ of $G$
satisfying $\dim G/H\leq 4n$ for $n\geq 22$, necessarily contain one
of the groups
\begin{align*}
\SO(n)\times \SO(4),~ \SO(n+3),~ \SO(n+2)\times \SO(2)
\end{align*}
as a factor---which includes into $\SO(n+4)$---due to table
\ref{PQKtable08}. (Note that whether $\SO(n+3)$ is a maximal rank
subgroup of $\SO(n+4)$ or not depends on the parity of $n$ being odd
or even.) By the same arguments we see that for $G=\SO(n+5)$ only
the following subgroups $H$ may appear:
\begin{align*}
&\SO(n+4),~ \SO(n+3)\times \SO(2),~ \SO(n+2)\times \SO(2),~
\SO(n+2)\times \SO(3)
\end{align*}
For $G=\SO(n+6)$ the group $H$ is out of the following list:
\begin{align*}
&\SO(n+5),~ \SO(n+4)\times \SO(2),~ \SO(n+3)\times \SO(3),~
\SO(n+3)\times \SO(2)
\end{align*}
In any of the cases there has to be an inclusion $\SO(n+k)$ for
$k\geq 0$ into $\Sp(n)$ by the isotropy
representation---cf.~\cite{kob1}, theorem VI.4.6.(2). By lemma
\ref{REClemma01}, however, this is impossible unless $k=0$. Thus we
see that
\begin{align*}
G/H=\frac{\SO(n+4)}{\SO(n)\times \SO(4)}=M
\end{align*}
since the action of $G$ thus necessarily is transitive.

\vspace{3mm}

Now suppose that $G=\Sp(\frac{n}{2}+2)$ for $n$ even. Again we use
table \ref{PQKtable08} to list maximal rank subgroups $H$ with $\dim
G/H \leq 4n$:
\begin{align*}
\Sp \bigg(\frac{n}{2}+1\bigg)\times \Sp(1),~ \Sp
\bigg(\frac{n}{2}+1\bigg)\times \U(1),
~\Sp\bigg(\frac{n}{2}\bigg)\times \Sp(2)
\end{align*}
If $H=\Sp(\frac{n}{2})\times \Sp(2)$, we see that $\dim G/H=4n$ and
that the action of $G$ is transitive. Thus $M$ is homogeneous and
symmetric. Yet, by the classification of Wolf spaces we derive that
\begin{align*}
M=\frac{\Sp\big(\frac{n}{2}+2\big)}{\Sp\big(\frac{n}{2}\big)\times
\Sp(2)}
\end{align*}
cannot be the case.

Now apply lemma \ref{REClemma02} in the other cases and derive that
the isotropy representation of $H$ is marked by a blockwise included
$\Sp(\frac{n}{2}+1)\hto{}\Sp(n)$ for $k>0$. This implies that the
sphere
\begin{align*}
\s^1\times \{1\} \times \stackrel{(n/2)}{\dots} \times \{1\} \hto{}
T^{n/2+1} \hto{} \Sp\bigg(\frac{n}{2}+1\bigg)
\end{align*}
is represented by $\s^1\times \{1\} \times \stackrel{(n-1)}{\dots}
\times \{1\}$ in the $\Sp(n)$-factor of the holonomy group
$\Sp(n)\Sp(1)$. Thus it fixes a quaternionic fixed-point component
of codimension $4$. Thus by theorem \cite{fang2}.1.2, p.~2, we
obtain that $M^{4n}\in\{\hh\pp^{4n},\Gr_2(\cc^{n+2})\}$.

If $G=\Sp(\frac{n+1}{2}+2)$ and $n$ is odd, virtually the same
arguments apply. That is, the list of isotropy subgroups $H$ is
given by
\begin{align*}
\Sp \bigg(\frac{n+1}{2}+1\bigg)\times \Sp(1),~ \Sp
\bigg(\frac{n+1}{2}+1\bigg)\times \U(1)
\end{align*}
As above this leads to a codimension four quaternionic
$\s^1$-fixed-point component.

\vspace{3mm}

Finally, suppose $n$ to be odd and the group $G$ to contain a factor
of the form $\Sp(\frac{n+1}{2}+1)$. The list of possible stabilisers
is given as
\begin{align*}
&\Sp\bigg(\frac{n+1}{2}\bigg)\times \Sp(1),~
\Sp\bigg(\frac{n+1}{2}\bigg)\times \U(1),~
\Sp\bigg(\frac{n+1}{2}-1\bigg)\times \Sp(2),
\\&\Sp\bigg(\frac{n+1}{2}-1\bigg)\times \Sp(1)\times \Sp(1)
\end{align*}
If $H=\Sp(\frac{n+1}{2}-1)\times \Sp(1)\times \Sp(1)$, the action of
$G$ again is transitive which is impossible by the classification of
Wolf spaces. If
\begin{align*}
H\in \bigg\{\Sp\bigg(\frac{n+1}{2}\bigg)\times
\Sp(1),~\Sp\bigg(\frac{n+1}{2}\bigg)\times \U(1)\bigg\}
\end{align*}
 we note that its $\Sp(\frac{n+1}{2})$-factor again maps into $\Sp(n)$ in a
blockwise way---cf.~lemma \ref{REClemma02}---by the isotropy
representation. This leads to a codimension four quaternionic
$\s^1$-fixed point component once more. Now suppose
\begin{align*}
H=\Sp\bigg(\frac{n+1}{2}-1\bigg)\times \Sp(2)
\end{align*}
Then the holonomy representation necessarily makes $H$ a subgroup of
$\Sp(n)$. By lemma \ref{REClemma04} this can only occur in the
standard blockwise way. Again this yields a quaternionic codimension
four $\s^1$-fixed-point component which leads to $M\in
\{\hh\pp^n,\Gr_2(\cc^{n+2})\}$.

\vspace{3mm}

In degree $n=21$ we see that similar arguments apply as for $n\geq
22$. However, we have that $\dim(\SO(25))=300<303=\dim \E_6\times
\F_4\times \Sp(1)$. Thus due to the assumption that $\dim
\Isom(M^{21})>303$ we may not identify the real Grassmannian
$\widetilde \Gr_4(\rr^{25})$ but only the quaternionic projective
space and the complex Grassmannian.
\end{proof}




\vfill

\begin{center}
\noindent
\begin{minipage}{\linewidth}
\small \noindent \textsc
{Manuel Amann} \\
\textsc{Fachbereich Mathematik und Informatik}\\
\textsc{Westf\"alische Wilhelms-Universit\"at M\"unster}\\
\textsc{Einsteinstr.~62}\\
\textsc{48149~M\"unster, Germany} \\[1ex]
\textsf{mamann@uni-muenster.de}\\
\textsf{http://wwwmath.uni-muenster.de/u/mamann}
\end{minipage}
\end{center}

\end{document}